\documentclass[10pt]{amsart} 
\usepackage{amscd,amssymb}
\usepackage{graphicx}
\usepackage[matrix,tips,arrow]{xy}
\newcommand{\choosegraphics}[2]{#2}  

\newtheorem{thm}[equation]{Theorem}
\newtheorem{cor}[equation]{Corollary}
\newtheorem{prop}[equation]{Proposition}
\newtheorem{lem}[equation]{Lemma}
\theoremstyle{definition}
\newtheorem{dfn}[equation]{Definition}
\newtheorem{rem}[equation]{Remark}

\newtheorem{notation}[equation]{Notation}
\newtheorem{prob}[equation]{Problem}
\numberwithin{equation}{section}

\newcommand{\iso}{\stackrel{\simeq}{\rightarrow}}
\newcommand{\inj}{\rightarrowtail}
\newcommand{\surj}{\twoheadrightarrow}
\newcommand{\opn}{\operatorname}
\newcommand{\cat}[1]{\operatorname{\mathsf{#1}}}

\newcommand{\rmitem}[1]{\item[\text{\textup{(#1)}}]}
\newcommand{\mfrak}[1]{\mathfrak{#1}}

\newcommand{\msf}[1]{\mathsf{#1}}
\newcommand{\mbf}[1]{\mathbf{#1}}
\newcommand{\mrm}[1]{\mathrm{#1}}
\newcommand{\mbb}[1]{\mathbb{#1}}

\newcommand{\tup}[1]{\textup{#1}}
\newcommand{\bsym}[1]{\boldsymbol{#1}}
\newcommand{\bra}[1]{\langle #1 \rangle}

\newcommand{\scrp}[1]{\scriptstyle{#1}}
\newcommand{\stick}{\rule[-1ex]{0ex}{4ex}}

\title[Derived Picard Groups]{Derived Picard Groups of Finite 
Dimensional Hereditary Algebras}
\author{Jun-ichi Miyachi and Amnon Yekutieli}
\date{19 July 2000}
\address{J. Miyachi: Department of Mathematics, Tokyo Gakugei
University, Koganei-shi, Tokyo, 184, Japan}
\email{miyachi@u-gakugei.ac.jp}
\address{A. Yekutieli: 
Department Mathematics and Computer
Science, Ben Gurion University, Be'er Sheva 84105, Israel}
\email{amyekut@math.bgu.ac.il}
\thanks{\textit{Keywords}: derived category; Picard group;
finite dimensional algebra; quiver.}
\subjclass{Primary: 18E30, 16G20; 
Secondary: 16G70, 16G60, 16E20, 14C22} 
\thanks{The second author was partially supported by the US-Israel 
Binational Science Foundation.}

\begin{document}

\begin{abstract}
Let $A$ be a finite dimensional algebra over a field $k$. The 
derived Picard group $\opn{DPic}_{k}(A)$ is the group of triangle 
auto-equivalences of $\cat{D}^{\mrm{b}}(\cat{mod} A)$ induced by 
two-sided tilting complexes.

We study the group $\opn{DPic}_{k}(A)$ when $A$ is hereditary 
and $k$ is algebraically closed. We obtain general results on the
structure of $\opn{DPic}_{k}(A)$, as well as explicit calculations for 
many cases, including all finite and tame representation types. 

Our method is to construct a representation of
$\opn{DPic}_{k}(A)$ on a certain infinite quiver
$\vec{\Gamma}^{\mrm{irr}}$. This representation is faithful when 
the quiver $\vec{\Delta}$ of $A$ is a tree, and then
$\opn{DPic}_{k}(A)$ is discrete. Otherwise a connected linear
algebraic group can occur as a factor of $\opn{DPic}_{k}(A)$.

When $A$ is hereditary, $\opn{DPic}_{k}(A)$ coincides with the 
full group of $k$-linear triangle auto-equivalences of 
$\cat{D}^{\mrm{b}}(\cat{mod} A)$. Hence we can calculate the group 
of such auto-equivalences for any triangulated category $\cat{D}$ 
equivalent to \linebreak
$\cat{D}^{\mrm{b}}(\cat{mod} A)$. These include the 
derived categories of piecewise hereditary algebras, and of 
certain noncommutative spaces introduced by Kontsevich-Rosen\-berg. 
\end{abstract}

\maketitle

\setcounter{section}{-1}
\section{Introduction and Statement of Results}

Let $k$ be a field and $A$ an associative unital $k$-algebra. We 
write $\cat{Mod} A$ for the category of left $A$-modules, and 
$\cat{D}^{\mrm{b}}(\cat{Mod} A)$ for the bounded derived category.
Let $A^{\circ}$ be the opposite 
algebra and $A^{\mrm{e}} := A \otimes_{k} A^{\circ}$
the enveloping algebra, so that 
$\cat{Mod} A^{\mrm{e}}$ is the category of $k$-central 
$A$-$A$-bimodules. 

A {\em two-sided tilting complex} a complex 
$T \in \cat{D}^{\mrm{b}}(\cat{Mod} A^{\mrm{e}})$
for which there exists another complex
$T^{\vee} \in \cat{D}^{\mrm{b}}(\cat{Mod} A^{\mrm{e}})$
satisfying
$T^{\vee} \otimes^{\mrm{L}}_{A} T \cong 
T \otimes^{\mrm{L}}_{A} T^{\vee} \cong A$.
This notion is due to Rickard \cite{Rd}.
The {\em derived Picard group} of $A$ \tup{(}relative to
$k$\tup{)} is
\[ \opn{DPic}_{k}(A) :=
\frac{ \{ \text{two-sided tilting complexes}\
T \in \cat{D}^{\mrm{b}}(\cat{Mod} A^{\mrm{e}}) \} }{
\text{isomorphism} } \]
with identity element $A$, product
$(T_{1}, T_{2}) \mapsto T_{1} \otimes_{A}^{\mrm{L}} T_{2}$
and inverse
$T \mapsto T^{\vee} := \opn{R} \opn{Hom}_{A}(T, A)$.
See \cite{Ye} for more details. 

Since every invertible bimodule is a two-sided tilting complex, 
$\opn{DPic}_{k}(A)$ contains the (noncommutative) Picard 
group $\opn{Pic}_{k}(A)$  as a subgroup. It also contains a 
central subgroup $\bra{\sigma} \cong \mbb{Z}$, where $\sigma$ is 
the class  of the two-sided tilting complex $A[1]$. In \cite{Ye} 
we showed that when $A$ is either local or commutative one has
$\opn{DPic}_{k}(A) = \opn{Pic}_{k}(A) \times \bra{\sigma}$.
This was discovered independently by Rouquier-Zimmermann 
\cite{Zi}, \cite{RZ}.
On the other hand in the smallest example of a $k$-algebra $A$
that is neither commutative nor local, namely the $2 \times 2$ upper
triangular matrix algebra, this equality fails. These observations
suggest that the group structure of $\opn{DPic}_{k}(A)$ should 
carry some information about the geometry of the noncommutative 
ring $A$.

This prediction is further motivated by another result in
\cite{Ye}, which says that $\opn{DPic}_{k}(A)$ classifies the
dualizing complexes over $A$. The geometric significance of 
dualizing complexes is well known (cf.\ \cite{RD} and \cite{YZ}).

 From a broader perspective, $\opn{DPic}_{k}(A)$ is related to the 
geometry of noncommutative schemes on the one hand, and to mirror 
symmetry and deformations of (commutative) smooth projective 
varieties on the other hand. See \cite{BO}, \cite{Ko}, \cite{KR} and 
\cite{Or}.

A good starting point for the study of the group $\opn{DPic}_{k}(A)$
is to consider {\em finite dimensional} $k$-algebras. The geometric 
object associated to a finite dimensional $k$-algebra $A$ is its 
quiver $\vec{\Delta}$, as defined by Gabriel (cf.\ \cite{GR} or
\cite{ARS}). It is worthwhile to note that from the point of view of 
noncommutative localization theory (cf.\ \cite{MR} Section 4.3) 
$\vec{\Delta}$ is the link graph of $A$. More on this in 
Remark \ref{rem1.1}. 

Some calculations of the groups $\opn{DPic}_{k}(A)$ for finite 
dimensional algebras have already been done. Let us mention the 
work of Rouquier-Zimmermann \cite{RZ} on Brauer tree algebras, and 
the work of Lenzing-Meltzer \cite{LM} on canonical algebras. 

In this paper we present a systematic study the group 
$\opn{DPic}_{k}(A)$ when $A$ is a finite dimensional {\em hereditary} 
algebra over an {\em algebraically closed} field $k$.
We obtain general results on the structure of $\opn{DPic}_{k}(A)$, 
as well as explicit calculations. These results carry over to piecewise 
hereditary algebras, as well as to certain noncommutative schemes. 
The rest of the Introduction is devoted to stating our main results.

The group $\opn{Aut}_{k}(A) = \opn{Aut}_{\cat{Alg} k}(A)$ 
of $k$-algebra automorphisms is a linear algebraic group over $k$, 
via the inclusion into 
$\opn{Aut}_{\cat{Mod} k}(A) = \opn{GL}(A)$.
This induces a structure of linear algebraic group on the quotient 
$\opn{Out}_{k}(A)$ of outer automorphisms. We denote by 
$\opn{Out}^{0}_{k}(A)$ the identity component of 
$\opn{Out}_{k}(A)$.

Recall that $A$ is a {\em basic} $k$-algebra if 
$A / \mfrak{r} \cong k \times \cdots \times k$, 
where $\mfrak{r}$ is the Jacobson radical. For a basic algebra one 
has $\opn{Out}_{k}(A) = \opn{Pic}_{k}(A)$.
A hereditary basic algebra $A$ is isomorphic to the {\em path 
algebra} $k \vec{\Delta}$ of its quiver. 
An algebra $A$ is {\em indecomposable} iff the quiver $\vec{\Delta}$ 
is connected.

For Morita equivalent $k$-algebras $A$ and $B$ one has 
$\opn{DPic}_{k}(A) \cong \opn{DPic}_{k}(B)$, and the quivers of $A$ 
and $B$ are isomorphic. According to a result of Brauer (see 
\cite{Po} Section 2) one has 
$\opn{Out}^{0}_{k}(A) \cong \opn{Out}^{0}_{k}(B)$. 
If $A \cong \prod_{i = 1}^{n} A_{i}$ then 
$\opn{DPic}_{k}(A) \cong
G \ltimes \prod_{i = 1}^{n} \opn{DPic}_{k}(A_{i})$,
where $G \subset S_{n}$ is a permutation group (cf.\ \cite{Ye} 
Lemma 2.6). Also 
$\vec{\Delta}(A) \cong \coprod \vec{\Delta}(A_{i})$
and
$\opn{Out}^{0}_{k}(A) \cong \prod \opn{Out}^{0}_{k}(A_{i})$.
Since the main result Theorem \ref{thm0.2} is stated 
in terms of $\vec{\Delta}$ and $\opn{Out}^{0}_{k}(A)$, 
we allow ourselves to assume throughout that $A$ is a basic 
indecomposable algebra. 

Given a quiver $\vec{Q}$ we denote by $\vec{Q}_{0}$ its vertex set. 
For a pair of vertices $x, y \in \vec{Q}_{0}$ we write $d(x, y)$ 
for the arrow-multiplicity, i.e.\ the number of arrows 
$\alpha: x \to y$. Let $\opn{Aut}(\vec{Q}_{0})$ be the permutation 
group of $\vec{Q}_{0}$, and let 
$\opn{Aut}(\vec{Q}_{0}; d)$ be the subgroup of arrow-multiplicity 
preserving permutations, namely
\[ \opn{Aut}(\vec{Q}_{0}; d) = \{ \pi \in \opn{Aut}(\vec{Q}_{0})
\mid d(\pi(x), \pi(y)) = d(x, y) \text{ for all } x, y \in 
\vec{Q}_{0}\} . \]

Write $\opn{Aut}(\vec{Q})$ for the automorphism group of the quiver 
$\vec{Q}$. Then $\opn{Aut}(\vec{Q}_{0}; d)$ is the image of the 
canonical homomorphism
$\opn{Aut}(\vec{Q}) \to \opn{Aut}(\vec{Q}_{0})$.
The surjection 
$\opn{Aut}(\vec{Q}) \surj \opn{Aut}(\vec{Q}_{0}; d)$
is split, and it is bijective iff 
$\vec{Q}$ has no multiple arrows. 

Of particular importance to us is a certain countable quiver 
$\vec{\Gamma}^{\mrm{irr}}$. 
This is a full subquiver of the Auslander-Reiten quiver 
$\vec{\Gamma}(\cat{D}^{\mrm{b}}(\cat{mod} A))$
of $\cat{D}^{\mrm{b}}(\cat{mod} A)$ as defined by Happel
\cite{Ha}. Here $\cat{mod} A$ is the category of finitely generated 
$A$-modules. If $A$ has finite representation type (i.e.\ 
$\vec{\Delta}$ is a Dynkin quiver) then 
$\vec{\Gamma}^{\mrm{irr}} \cong \vec{\mbb{Z}} \vec{\Delta}$, 
where $\vec{\mbb{Z}} \vec{\Delta}$ is the quiver introduced by 
Riedtmann \cite{Rn}. Otherwise
$\vec{\Gamma}^{\mrm{irr}} \cong \mbb{Z} \times \vec{\mbb{Z}} 
\vec{\Delta}$. See Definitions \ref{dfn2.1} and \ref{dfn2.2}
for the definition of the quivers
$\vec{\Gamma}^{\mrm{irr}}$ and $\vec{\mbb{Z}} \vec{\Delta}$, and 
see Figures \ref{fig2} and \ref{fig3} for illustrations. The group 
$\opn{DPic}_{k}(A)$ acts on $\vec{\Gamma}^{\mrm{irr}}_{0}$ by 
arrow-multiplicity preserving permutations, giving rise to a group 
homomorphism
$q: \opn{DPic}_{k}(A) \to 
\opn{Aut}(\vec{\Gamma}^{\mrm{irr}}_{0}; d)$.

Define the bimodule $A^{*} := \opn{Hom}_{k}(A, k)$. Then $A^{*}$ is 
a two-sided tilting complex, the functor 
$M \mapsto A^{*} \otimes^{\mrm{L}}_{A} M \cong
\mrm{R} \opn{Hom}_{A}(M, A^{*})$ 
is the {\em Serre functor} of 
$\cat{D}^{\mrm{b}}(\cat{mod} A)$ in the sense of \cite{BK}, and 
$M \mapsto A^{*}[-1] \otimes^{\mrm{L}}_{A} M$
is the {\em translation functor} in the sense of \cite{Ha} Section
I.4. We write $\tau \in \opn{DPic}_{k}(A)$ for the element 
represented by $A^{*}[-1]$. Then $\tau$ is translation of the quiver 
$\vec{\Gamma}^{\mrm{irr}}$. Let us denote by
$\opn{Aut}(\vec{\Gamma}^{\mrm{irr}}_{0}; d)^{\bra{\tau, \sigma}}$
the subgroup of 
$\opn{Aut}(\vec{\Gamma}^{\mrm{irr}}_{0}; d)$
consisting of permutations that commute with $\tau$ and $\sigma$.

Here is the main result of the paper.

\begin{thm} \label{thm0.2}
Let $A$ be an indecomposable basic hereditary finite dimensional 
algebra over an algebraically closed field $k$, with quiver 
$\vec{\Delta}$.
\begin{enumerate}
\item There is an exact sequence of groups
\[ 1 \to \opn{Out}^{0}_{k}(A) \to \opn{DPic}_{k}(A) 
\xrightarrow{q} 
\opn{Aut}(\vec{\Gamma}^{\mrm{irr}}_{0} ; d)^{\bra{\tau, \sigma}}
\to 1 . \]
This sequence splits.
\item If $A$ has finite representation type then there is an 
isomorphism of groups
\[ \opn{DPic}_{k}(A) \cong
\opn{Aut}(\vec{\mbb{Z}} \vec{\Delta})^{\bra{\tau}} . \]
\item If $A$ has infinite representation type then 
there is an isomorphism of groups
\[ \opn{DPic}_{k}(A) \cong
\bigl( \opn{Aut}((\vec{\mbb{Z}} \vec{\Delta})_{0}; d)^{\bra{\tau}}
\ltimes \opn{Out}^{0}_{k}(A) \bigr) \times \mbb{Z} . \]
\end{enumerate}
\end{thm}

The factor $\mbb{Z}$ of $\opn{DPic}_{k}(A)$ in part 3 is 
generated by $\sigma$. If $\vec{\Delta}$ has no multiple arrows 
then so does $\vec{\mbb{Z}} \vec{\Delta}$, and hence
$\opn{Aut}((\vec{\mbb{Z}} \vec{\Delta})_{0}; d) = 
\opn{Aut}(\vec{\mbb{Z}} \vec{\Delta})$.
The proof of Theorem \ref{thm0.2} is in Section 3
where it is stated again as Theorem \ref{thm3.2}.

Recall that a finite dimensional $k$-algebra $B$ is called 
{\em piecewise hereditary} of type $\vec{\Delta}$ if 
$\cat{D}^{\mrm{b}}(\cat{mod} B) \approx 
\cat{D}^{\mrm{b}}(\cat{mod} A)$
where $A = k \vec{\Delta}$ for some finite quiver $\vec{\Delta}$ 
without oriented cycles. By \cite{Rd} Corollary 3.5 one knows that 
$\opn{DPic}_{k}(B) \cong \opn{DPic}_{k}(A)$. The next corollary 
follows. 

\begin{cor} \label{cor0.2}
Suppose $B$ is a piecewise hereditary $k$-algebra of type 
$\vec{\Delta}$. Then $\opn{DPic}_{k}(B)$ is
described by Theorem \tup{\ref{thm0.2}} with $A = k \vec{\Delta}$.
\end{cor}

In Section 4 we work out explicit descriptions of the groups
$\opn{Pic}_{k}(A)$ and $\opn{DPic}_{k}(A)$ for the Dynkin and affine
quivers, as well as for some wild quivers with multiple arrows. 
As an example we present below the explicit description of  
$\opn{DPic}_{k}(A)$ for a Dynkin quiver of type $A_{n}$ 
(which corresponds to upper triangular $n \times n$ matrices). The 
corollary is extracted from Theorem \ref{thm4.1}. 

\begin{cor} \label{cor0.4}
Suppose $\vec{\Delta}$ is a Dynkin quiver of type $A_{n}$ and 
$A = k \vec{\Delta}$. Then
$\opn{DPic}_{k}(A)$ is an abelian group generated by
$\tau$ and $\sigma$, with one relation
\[ \tau^{n + 1} = \sigma^{-2} . \]
\end{cor}

The relation $\tau^{n + 1} = \sigma^{-2}$ was already discovered 
by E. Kreines (cf.\ \cite{Ye} Appendix).
This relation has been known also to Kontsevich, and in 
his terminology $\cat{D}^{\mrm{b}}(\cat{mod} A)$
is ``fractionally Calabi-Yau of dimension $\frac{n-1}{n+1}$''
(see \cite{Ko}; note that the Serre functor is $\tau \sigma$).

Suppose $\cat{D}$ is a $k$-linear triangulated category that's 
equivalent to a small category. Denote by 
$\opn{Out}_{k}^{\mrm{tr}}(\cat{D})$ the group of $k$-linear 
triangle auto-equivalences of $\cat{D}$ modulo functorial 
isomorphisms. For a finite dimensional algebra $A$ one has 
$\opn{DPic}_{k}(A) \subset 
\opn{Out}_{k}^{\mrm{tr}}(\cat{D}^{\mrm{b}}(\cat{mod} A)) $, 
with equality when $A$ is hereditary (cf.\ Corollary 
\ref{cor1.5}). 

In \cite{KR} Kontsevich-Rosenberg introduce the noncommutative 
projective space $\mbf{NP}^{n}_{k}$, $n \geq 1$. They state that 
$\cat{D}^{\mrm{b}}(\cat{Coh} \mbf{NP}^{n}_{k})$ 
is equivalent to 
$\cat{D}^{\mrm{b}}(\cat{mod} k \vec{\Omega}_{n+1})$,
where $\vec{\Omega}_{n+1}$ is the quiver in Figure \ref{fig5}, 
and $\cat{Coh} \mbf{NP}^{n}_{k}$ is the category of 
coherent sheaves. 
By Beilinson's results in \cite{Be}, there is an equivalence
$\cat{D}^{\mrm{b}}(\cat{Coh} \mbf{P}^{1}_{k}) \approx
\cat{D}^{\mrm{b}}(\cat{mod} k \vec{\Omega}_{2})$. 
Combining Theorem \ref{thm4.3} and Corollary \ref{cor1.5} we 
get the next corollary.

\begin{cor}
Let $X$ be either $\mbf{NP}^{n}_{k}$ \tup{(}$n \geq 1$\tup{)} or 
$\mbf{P}^{n}_{k}$ \tup{(}$n = 1$\tup{)}. Then
\[ \opn{Out}_{k}^{\mrm{tr}}(\cat{D}^{\mrm{b}}(\cat{Coh} X))
\cong 
\mbb{Z} \times \bigl( \mbb{Z} \ltimes \opn{PGL}_{n+1}(k) \bigr) . \]
\end{cor}
 
In Section 5 we look at a tree $\Delta$ with $n$ 
vertices. Every orientation $\omega$ of $\Delta$ gives a quiver 
$\vec{\Delta}_{\omega}$. The equivalences between the various 
categories 
$\cat{D}^{\mrm{b}}(\cat{mod} k \vec{\Delta}_{\omega})$
form the derived Picard groupoid 
$\opn{DPic}_{k}(\Delta)$. The subgroupoid generated by the 
two-sided tilting complexes of \cite{APR} is called the 
{\em reflection groupoid} $\opn{Ref}(\Delta)$. We show that there 
is a surjection $\opn{Ref}(\Delta) \surj W(\Delta)$, 
where $W(\Delta) \subset \opn{GL}_{n}(\mbb{Z})$ is the Weyl group 
as in \cite{BGP}. We also prove that for any orientation $\omega$, 
$\opn{Ref}(\Delta)(\omega,\omega) = \bra{\tau_{\omega}}$
where $\tau_{\omega} \in \opn{DPic}_{k}(A_{\omega})$ 
is the translation. 
 
\medskip \noindent \textbf{Acknowledgments.}\
We wish to thank A. Bondal, I. Reiten and M. Van den Bergh for 
very helpful conversations and correspondences. Thanks to the 
referee for suggestions and improvements to the paper. 
Some of the work on the paper was done during visits to MIT and 
the University of Washington, and we thank the departments of 
mathematics at these universities for their hospitality. The 
second author was supported by the Weizmann Institute of Science 
throughout most of this research.

\section{Conventions and Preliminary Results}

In this section we fix notations and conventions to be used 
throughout the paper. This is needed since there are conflicting 
conventions in the literature regarding quivers and path algebras. 
We also prove two preliminary results. 

Throughout the paper $k$ denotes a fixed algebraically closed field.
Our notation for a quiver is $\vec{Q} = (\vec{Q}_{0}, \vec{Q}_{1})$;
$\vec{Q}_{0}$ is the set of vertices, and $\vec{Q}_{1}$ is the set 
of arrows. For $x, y \in \vec{Q}_{0}$, $d(x, y)$ denotes the 
number of arrows $x \to y$. 

In this section the letter $\msf{A}$ denotes a $k$-linear category 
that's equivalent to a small full subcategory of itself (this 
assumption avoids some set theoretical problems). 
Let us write $\opn{Aut}_{k}(\msf{A})$ for the class of $k$-linear 
auto-equivalences of $\msf{A}$. Then the set
\begin{equation} \label{eqn1.4}
\opn{Out}_{k}(\msf{A}) = \frac{\opn{Aut}_{k}(\msf{A})}{ 
\text{functorial isomorphism}}
\end{equation}
is a group.

Suppose $\msf{A}$ is a $k$-linear additive Krull-Schmidt category
(i.e.\ $\opn{dim}_{k} \opn{Hom}_{\msf{A}}(M, N)$ \linebreak 
$< \infty$
and all idempotents split). 
We define the quiver $\vec{\Gamma}(\msf{A})$
of $\msf{A}$ as follows: $\vec{\Gamma}_{0}(\msf{A})$ is the 
set of isomorphism classes of indecomposable objects of $\msf{A}$. 
For two vertices $x, y$ there are $d(x, y)$ arrows 
$\alpha : x \to y$, where we choose representatives $M_{x} \in x$, 
$M_{y} \in y$, 
$\opn{Irr}(M_{x}, M_{y}) = \opn{rad}(M_{x}, M_{y}) /
\opn{rad}^{2}(M_{x}, M_{y})$ 
is the space of irreducible morphisms
and $d(x, y) := \opn{dim}_{k} \opn{Irr}(M_{x}, M_{y})$.
See \cite{Rl} Section 2.2 for full details.

If $\msf{A}$ is a $k$-linear category (possibly without direct sums)
we can embed it in the additive category
$\msf{A} \times \mbb{N}$, where a morphism
$(x, m) \to (y, n)$ is an $n \times m$ matrix with entries in
$\msf{A}(x, y) = \opn{Hom}_{\msf{A}}(x, y)$.
Of course if $\msf{A}$ is additive then
$\msf{A} \approx \msf{A} \times \mbb{N}$.
If $\msf{A} \times \mbb{N}$ is Krull-Schmidt then we shall write
$\vec{\Gamma}(\msf{A})$ for the quiver
$\vec{\Gamma}(\msf{A} \times \mbb{N})$.

Let $\vec{Q}$ be a quiver. Assume that for every vertex
$x \in \vec{Q}_{0}$ the number of arrows starting or ending at 
$x$ is finite, and for every two vertices $x, y \in \vec{Q}_{0}$ 
there is only a finite number of oriented paths from $x$ to $y$.
Let $k \bra{\vec{Q}}$ be the {\em path category}, whose set of 
objects is $\vec{Q}_{0}$, the morphisms are generated by the 
identities and the arrows, and the only relations arise from 
incomposability of paths. Observe that this differs from the 
definition in \cite{Rl}, where the path category corresponds to 
$k \bra{\vec{Q}} \times \mbb{N}$ in our notation. 
The morphism spaces of $k \bra{\vec{Q}}$ are $\mbb{Z}$-graded, 
where the arrows have degree $1$. If $I \subset k \bra{\vec{Q}}$ 
is any ideal contained in 
$\opn{rad}^{2}_{k \bra{\vec{Q}}} = \bigoplus _{n \geq 2} 
k \bra{\vec{Q}}_{n}$, 
and $k \bra{\vec{Q}, I} := k \bra{\vec{Q}} / I$ 
is the quotient category, then the 
additive category $k \bra{\vec{Q}, I} \times \mbb{N}$ is 
Krull-Schmidt, and the quiver of $k \bra{\vec{Q}, I}$ is
$\vec{\Gamma}(k \bra{\vec{Q}, I}) = \vec{Q}$.

Let $A$ be a finite dimensional $k$-algebra. In representation 
theory there are three equivalent ways to define the quiver 
$\vec{\Delta} = \vec{\Delta}(A)$ of $A$. The set 
$\vec{\Delta}_{0}$ enumerates either a complete set of primitive 
orthogonal idempotents 
$\{ e_{x} \}_{x \in \vec{\Delta}_{0}}$, as in \cite{ARS} Section 
III.1; or it enumerates the simple $A$-modules 
$\{ S_{x} \}_{x \in \vec{\Delta}_{0}}$, as in \cite{Rl} Section 
2.1; or it enumerates the indecomposable projective $A$-modules 
$\{ P_{x} \}_{x \in \vec{\Delta}_{0}}$, as in \cite{Rl} Section 
2.4. The arrow multiplicity is in all cases
\[ d(x, y) = 
\opn{dim}_{k} e_{x} (\mfrak{r} / \mfrak{r}^{2}) e_{y} = 
\opn{dim}_{k} \opn{Ext}^{1}_{A}(S_{y}, S_{x}) = 
\opn{dim}_{k} \opn{Irr}_{\cat{proj} A}(P_{x}, P_{y}) . \]
Here $\mfrak{r}$ is the Jacobson radical and $\cat{proj} A$ is the 
category of finitely generated projective modules, which is 
Krull-Schmidt. Observe that the third definition is just
$\vec{\Delta}(A) = \vec{\Gamma}(\cat{proj} A)$. 

\begin{rem} \label{rem1.1}
The set $\vec{\Delta}_{0}$ also enumerates the prime spectrum of $A$, 
$\opn{Spec} A \cong \{ \mfrak{p}_{x} \}_{x \in \vec{\Delta}_{0}}$. 
One can show that 
$\mfrak{r} / \mfrak{r}^{2} \cong 
\bigoplus_{x, y \in \vec{\Delta}_{0}}
(\mfrak{p}_{x} \cap \mfrak{p}_{y}) / \mfrak{p}_{x} \mfrak{p}_{y}$
as $A$-$A$-bimodules. This implies that $d(x, y) > 0$ iff
there is a second layer link 
$\mfrak{p}_{x} \rightsquigarrow \mfrak{p}_{y}$
(cf.\ \cite{MR} Section 4.3.7). Thus if we ignore multiple arrows, 
the quiver $\vec{\Delta}$ is precisely the link graph of $A$. 
\end{rem}

Recall that a {\em translation} $\tau$ is an injective function from 
a subset of $\vec{Q}_{0}$, called the set of non-projective vertices, 
to $\vec{Q}_{0}$, such that $d(\tau(y), x) = d(x, y)$. 
$\vec{Q}$ is a stable translation quiver if it comes with a 
translation $\tau$ such that all vertices are non-projective.
A {\em polarization} $\mu$ is an injective function defined on the 
set of arrows $\beta: x \to y$ ending in non-projective vertices, 
with $\mu(\beta) : \tau(y) \to x$. Cf.\ \cite{Rl} Section 2.2.

\begin{notation} \label{not1.1}
Suppose the quiver $\vec{Q}$ has a translation $\tau$ and a 
polarization $\mu$. Given a non-projective vertex $y \in \vec{Q}_{0}$ 
let $x_{1}, \ldots, x_{m}$ be some labeling, without repetition, 
of the set of vertices 
$\{ x \mid \text{there is an arrow } x \to y \}$ . 
Correspondingly label the arrows
$\beta_{i, j} : x_{i} \to y$
and
$\alpha_{i, j} : \tau(y) \to x_{i}$,
where $i = 1, \ldots, m$;
$j = 1, \ldots, d_{i} = d(x_{i}, y)$; and
$\alpha_{i, j} = \mu(\beta_{i, j})$.
The {\em mesh ending at} $y$ is the subquiver with vertices 
$\{ \tau(y), x_{i}, y \}$ and arrows 
$\{ \alpha_{i, j}, \beta_{i, j} \}$.
\end{notation}

\begin{figure}
\choosegraphics{
\[ \UseTips \begin{xy} 
            (0,0)*+@{*}="ty"*+!CR{\tau(y)},"ty"
\ar@{->}^{\alpha_{1}}    (20,15)*+@{*}="x1"*+!DL{x_{1}},"x1"
\ar@{}      (20,7.5)*+@{}*!C{\vdots}
\ar@{->}^{\alpha_{i}}    (20,0)*+@{*}="xi"*+!DL{x_{i}},"xi"
\ar@{}      (20,-7.5)*+@{}*!C{\vdots}
\ar@{->}^{\alpha_{m}}    (20,-15)*+@{*}="xm"*+!UL{x_{m}},"xm"
\ar@{->}^{\beta_{1}}   "x1";(40,0)*++@{*}="y"*+!CL{y},"y"
\ar@{->}^{\beta_{i}}   "xi";"y"
\ar@{->}^{\beta_{m}}   "xm";"y"
\end{xy} \]
}{
\includegraphics[clip]{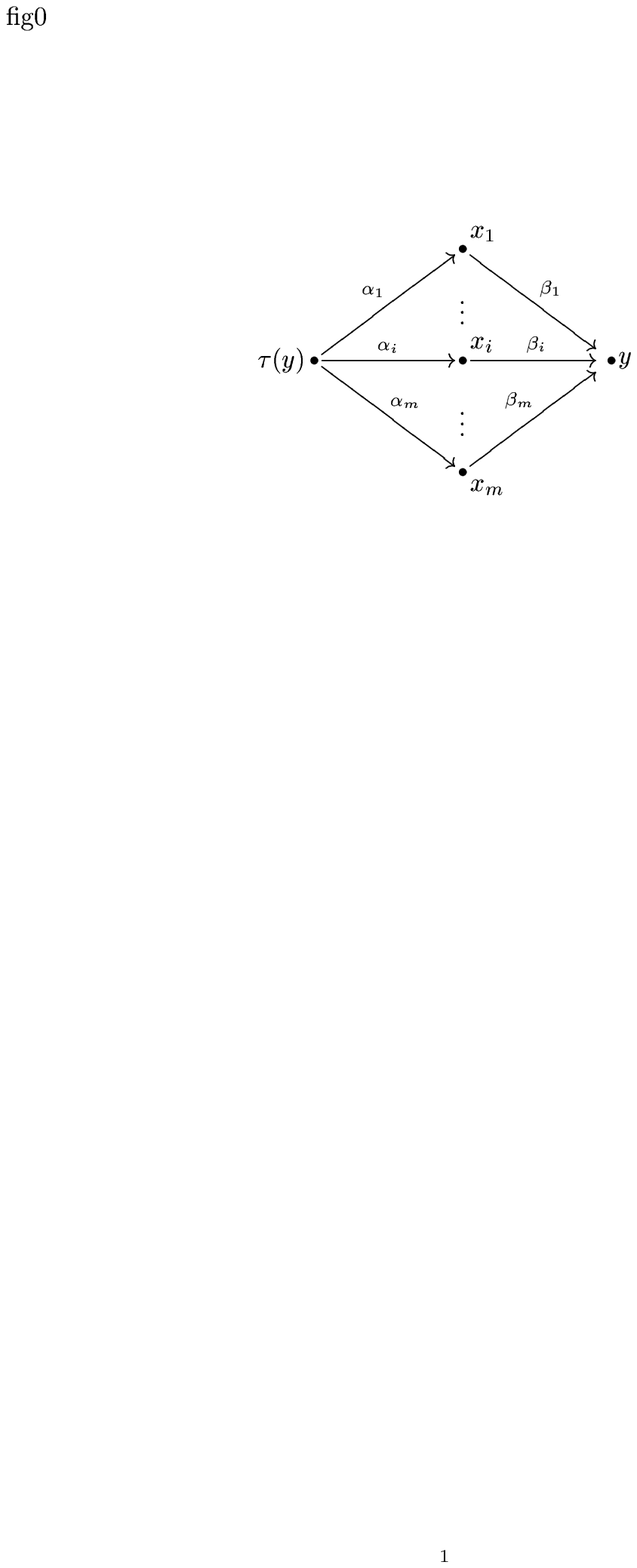}
}
\caption{The mesh ending at the vertex $y$ when $d_{i}=1$.} 
\label{fig0}
\end{figure}

If $\vec{Q}$ has no multiple arrows then $d_{i} = 1$ and the 
picture of the mesh ending at $y$ is shown in Figure \ref{fig0}.

The {\em mesh relation} at $y$ is defined to be
\begin{equation} \label{eqn1.1}
\sum_{i = 1}^{m} \sum_{j = 1}^{d_{i}}
\beta_{i, j} \alpha_{i, j} \in  
\opn{Hom}_{k \bra{\vec{Q}}} \bigl( \tau(y), y \bigr) . 
\end{equation}
It is a homogeneous morphism of degree $2$.

\begin{dfn} \label{dfn1.1}
Let $I_{\mrm{m}}$ be the mesh ideal in the category 
$k \bra{\vec{Q}}$,
i.e.\ the two sided ideal generated by the mesh relations 
(\ref{eqn1.1}) where $y$ runs over all non-projective vertices. 
The quotient category 
\[ k \bra{\vec{Q}, I_{\mrm{m}}} := k \bra{\vec{Q}} / I_{\mrm{m}}
\]
is called the {\em mesh category}. 
\end{dfn}

Observe that in \cite{Rn}, \cite{Rl} and \cite{Ha}
the notation for $k \bra{\vec{Q}, I_{\mrm{m}}}$ is 
$k(\vec{Q})$.

Now let $\vec{\Delta}$ be a finite quiver without oriented cycles and 
$A = k \vec{\Delta}$ the path algebra. Our convention for the 
multiplication in $A$ is as follows. If
$x \xrightarrow{\alpha} y$ and $y \xrightarrow{\beta} z$ are
paths in $\vec{\Delta}$, and if
$x \xrightarrow{\gamma} z$ is the concatenated path, then
$\gamma = \alpha \beta$ in $A$.
We note that the composition rule in the path category 
$k \bra{\vec{\Delta}}$ is opposite to that in $A$, so that 
$\bigoplus_{x, y} \opn{Hom}_{k \bra{\vec{\Delta}}}(x, y) =
A^{\circ}$.

For every $x \in \vec{\Delta}_{0}$ let $e_{x} \in A$ be the
corresponding idempotent, and let $P_{x} = A e_{x}$ be the
indecomposable projective $A$-module. So
$\{ P_{x} \}_{x \in \vec{\Delta}_{0}}$
is a set of representatives of the isomorphism classes of
indecomposable projective $A$-modules. Define
$\msf{P} \subset \cat{mod} A$ to be the full subcategory on the
objects $\{ P_{x} \}_{x \in \vec{\Delta}_{0}}$.
Then
$\msf{P} \times \mbb{N} \approx \cat{proj} A$
and
$\vec{\Delta} \cong \vec{\Gamma}(\msf{P}) \cong
\vec{\Gamma}(\cat{proj} A)$.

There is an equivalence of categories
$k \bra{\vec{\Delta}} \xrightarrow{\approx} \msf{P}$
that sends $x \mapsto P_{x}$, and an arrow
$\alpha : x \to y$ goes to the right multiplication
$P_{x} = A e_{x} \xrightarrow{\alpha} P_{y} = A e_{y}$. 
We will identify $\msf{P}$ and $k \bra{\vec{\Delta}}$ in this way. 

Recall that the automorphism group $\opn{Aut}_{k}(A)$ is a linear 
algebraic group. Let $H$ be the closed subgroup
\[ H := \{ F \in \opn{Aut}_{k}(A) \mid F(e_{x}) = e_{x}
\text{ for all } x \in \vec{\Delta}_{0} \} . \]

\begin{lem} \label{lem1.1}
$H$ is connected.
\end{lem}

\begin{proof}
For each pair $x, y \in \vec{\Delta}_{0}$ the $k$-vector space
$\msf{P}(x,y) := \opn{Hom}_{\msf{P}}(x, y) \cong 
e_{x} A e_{y}$ 
is graded. Let $\msf{P}(x,y)_{i}$ be the homogeneous component of 
degree $i$, and
\[ Y := \prod_{x, y \in \vec{\Delta}_{0}}
\Bigl( \opn{Aut}_{k} 
\bigl( \msf{P}(x, y)_{1} \bigr) \times
\opn{Hom}_{k} \bigl( \msf{P}(x, y)_{1}, 
\msf{P}(x, y)_{\geq 2} \bigr) \Bigr) . \]
This is a connected algebraic variety. Since $A$ is generated as 
$k$-algebra by the idempotents and the arrows, and the only relations 
in $A$ are the monomial relations arising from incomposability of 
paths, it follows that any element $F' \in Y$ extends 
uniquely to a $k$-algebra automorphism $F$ of $A$ that fixes the 
idempotents. Conversely any automorphism $F \in H$ restricts to an 
element $F'$ of $Y$. This bijection $Y \to H$ is an isomorphism of 
varieties. Hence $H$ is connected.
\end{proof}

The next result is partially proved in \cite{GS} Theorem 4.8 (they 
assume $k$ has characteristic $0$). 

\begin{prop} \label{prop1.1}
Let $A$ be a basic hereditary finite dimensional algebra over an 
algebraically closed field $k$, with quiver $\vec{\Delta}$.
\begin{enumerate}
\item There is a split exact sequence of groups
\[ 1 \to \opn{Out}^{0}_{k}(A) \to \opn{Pic}_{k}(A) \to
\opn{Aut}(\vec{\Delta}_{0};d) \to 1 . \]
\item The group $\opn{Out}^{0}_{k}(A)$ is trivial when 
$\vec{\Delta}$ is a tree.
\end{enumerate}
\end{prop}

\begin{proof}
1. Since $A$ is basic we have 
$\opn{Out}_{k}(A) = \opn{Pic}_{k}(A)$.
By Morita theory we have
$\opn{Pic}_{k}(A) \cong \opn{Out}_{k}(\cat{Mod} A)$. Any 
auto-equivalence of the category $\msf{P}$ extends to an 
auto-equivalence of $\cat{Mod} A$ (using projective resolutions), 
and this induces an isomorphism of groups
$\opn{Out}_{k}(\msf{P}) \iso \opn{Out}_{k}(\cat{Mod} A)$. 

The class of auto-equivalences $\opn{Aut}_{k}(\msf{P})$ is actually 
a group here. In fact $\opn{Aut}_{k}(\msf{P})$ can be identified with 
the subgroup of $\opn{Aut}_{k}(A)$ consisting of automorphisms that 
permute the set of idempotents $\{ e_{x} \} \subset A$.

Define a homomorphism of groups
$q: \opn{Out}_{k}(A) \to \opn{Aut}(\vec{\Delta}_{0};d)$
by $q(F)(x) = y$ if $F P_{x} \cong P_{y}$. Thus we get a commutative 
diagram
\[ \UseTips \xymatrix@M+1ex{
H \ar@{>->}[r] 
& \opn{Aut}_{k}(\msf{P}) \ar@{>->}[r] \ar@{->>}[d]^{f}
& \opn{Aut}_{k}(A) \ar@{->>}[d]^{g} \\
& \opn{Out}_{k}(\msf{P}) \ar[r]^{\cong} 
& \opn{Out}_{k}(A) \ar[r]^{q} &
\opn{Aut}(\vec{\Delta}_{0};d) . } \]
For an element $F \in \opn{Aut}_{k}(\msf{P})$ we have
$F|_{ \{e_{x}\} } = q f(F)$, and hence 
$\opn{Ker}(q) = f(H) = g(H)$. 
According to Lemma \ref{lem1.1}, $H$ is connected. Because $g$ is a 
morphism of varieties we see that $\opn{Ker}(q)$ is connected. 
But the index of $\opn{Ker}(q)$ is finite, so we get 
$\opn{Ker}(q) = \opn{Out}^{0}_{k}(A)$.

In order to split $q$ we choose any splitting of 
$\opn{Aut}(\vec{\Delta}) \surj \opn{Aut}(\vec{\Delta}_{0};d)$
and compose it with the homomorphism
$\opn{Aut}(\vec{\Delta}) \to \opn{Aut}_{k}(\msf{P})$.

\medskip \noindent 2.
When $\vec{\Delta}$ is a tree the group $H$ is a torus:
$H \cong \prod_{x, y \in \vec{\Delta}_{0}}
\opn{Aut}_{k} \bigl( \msf{P}(x, y)_{1} \bigr)$.
In fact $H$ consists entirely of inner automorphisms that are
conjugations by elements of the form
$\sum \lambda_{x} e_{x}$ with $\lambda_{x} \in k^{\times}$. Thus 
$g(H) = 1$.
\end{proof}

The next theorem seems to be known to some experts, but we could not 
locate any reference in the literature. Since it is needed in the 
paper we have included a short proof. For a left coherent ring $A$ 
(e.g.\ a hereditary ring) we denote by $\cat{mod} A$ the category 
of coherent $A$-modules. In the theorem $k$ could be any field. 

\begin{thm} \label{thm1.4}
Suppose $A$ is a hereditary $k$-algebra. Then any $k$-linear 
triangle auto-equivalence of $\cat{D}^{\mrm{b}}(\cat{mod} A)$
is standard.
\end{thm}

\begin{proof}
Let $F$ be a $k$-linear triangle auto-equivalence of 
$\cat{D}^{\mrm{b}}(\cat{mod} A)$. By \cite{Rd} Corollary 3.5
there exits a two-sided tilting complex $T$ with 
$T \cong F A$ in $\cat{D}^{\mrm{b}}(\cat{mod} A)$. Replacing $F$ 
with $(T^{\vee} \otimes^{\mrm{L}}_{A} -) F$ 
we may assume that $F A \cong A$. Hence
$F(\cat{mod} A) \subset \cat{mod} A$, and $F|_{\cat{mod} A}$ is an 
equivalence. Classical Morita theory says that 
$F|_{\cat{mod} A} \cong (P \otimes_{A} -)$ for some invertible 
bimodule $P$. So replacing $F$ by $(P^{\vee} \otimes_{A} -) F$ we can 
assume that there is an isomorphism
$\phi^{0}: F|_{\cat{mod} A} \cong \mbf{1}_{\cat{mod} A}$.

Now for every object $M \in \cat{D}^{\mrm{b}}(\cat{mod} A)$ 
we can choose an isomorphism 
$M \cong \bigoplus_{i} M_{i}[-i]$ with $M_{i} \in \cat{mod} A$
(cf.\ \cite{Ha} Lemma I.5.2). 
Define $\phi_{M}: F M \iso M$ to be the composition
\[ F M \cong \bigoplus (F M_{i})[-i] 
\xrightarrow{\sum \phi^{0}_{M_{i}}[-i]} 
\bigoplus_{i} M_{i}[-i] \cong M . \]
According to the proof of \cite{BO} Proposition A.3, for any 
morphism $\alpha: M \to N$ one has
$\phi_{N} F(\alpha) = \alpha \phi_{M}$, so
$\phi: F \to \mbf{1}_{\cat{D}^{\mrm{b}}(\cat{mod} A)}$ 
is an isomorphism of functors. 
\end{proof}

\begin{cor} \label{cor1.5}
Suppose $A$ is a hereditary $k$-algebra. Then
\[ \opn{DPic}_{k}(A) \cong 
\opn{Out}_{k}^{\mrm{tr}}(\cat{D}^{\mrm{b}}(\cat{mod} A)) . \]
\end{cor}

\begin{proof}
The group homomorphism
$\opn{DPic}_{k}(A) \to 
\opn{Out}_{k}^{\mrm{tr}}(\cat{D}^{\mrm{b}}(\cat{mod} A))$
is injective, say by \cite{Ye} Proposition 2.2, and it is 
surjective by the theorem.
\end{proof}

\section{An Equivalence of Categories}

In this section we prove the technical result Theorem \ref{thm2.3}.
It holds for any finite dimensional hereditary $k$-algebra $A$. In 
the special case of finite representation type, Theorem \ref{thm2.3}
is just \cite{Ha} Proposition I.5.6. Our result is the derived 
category counterpart of \cite{Rl} Lemma 2.3.3. For notation see 
Section 1 above.

We use a few facts about Auslander-Reiten triangles in  
$\cat{D}^{\mrm{b}}(\cat{mod} A)$. These facts are well known to 
experts in representation theory, but for the benefit of other 
readers we have collected them in Theorems \ref{thm2.1} and 
\ref{thm2.2}.

Let $\cat{D}$ be a $k$-linear triangulated category, which is 
Krull-Schmidt (as additive category).
As in any Krull-Schmidt category, sink and source morphisms
can be defined in $\cat{D}$; cf.\ \cite{Rl} Section 2.2. 
In \cite{Ha} Section I.4, 
Happel defines Auslander-Reiten triangles in $\cat{D}$, 
generalizing the Aus\-lander-Reiten (or almost split) sequences 
in an abelian Krull-Schmidt category. A triangle
$M' \xrightarrow{g} M \xrightarrow{f} M'' \to M'[1]$ 
in $\cat{D}$ is an Auslander-Reiten triangle if $g$ is a source 
morphism, or equivalently if $f$ is a sink morphism. 
As before, we denote by $M_{x} \in \cat{D}$ an indecomposable object
in the isomorphism class $x \in \vec{\Gamma}(\cat{D})$.

Now let  $\vec{\Delta}$ be a finite quiver without oriented 
cycles, and $A = k \vec{\Delta}$ the path algebra.
For $M \in \cat{mod} k$ let 
$M^{*} :=  \opn{Hom}_{k}(M, k)$. 
Define auto-equivalences $\sigma$ and $\tau$ of 
$\cat{D}^{\mrm{b}}(\cat{mod} A)$
by $\sigma M := M[1]$
and
$\tau M := \opn{R} \opn{Hom}_{A}(M, A)^{*}[-1] \cong
A^{*}[-1] \otimes^{\mrm{L}}_{A} M$.

\begin{thm}[Happel, Ringel] \label{thm2.1}
Let $A = k \vec{\Delta}$. Then the following hold.
\begin{enumerate}
\item As an additive $k$-linear category, 
$\cat{D}^{\mrm{b}}(\cat{mod} A)$ is a Krull-Schmidt category.
\item The quiver 
$\vec{\Gamma} := \vec{\Gamma}(\cat{D}^{\mrm{b}}(\cat{mod} A))$
is a stable translation quiver, and the translation $\tau$ 
satisfies $M_{\tau(x)} \cong \tau M_{x}$.
\item The Auslander-Reiten triangles in 
$\cat{D}^{\mrm{b}}(\cat{mod} A)$ \tup{(}up to isomorphism\tup{)}
correspond bijectively  to the meshes in $\vec{\Gamma}$. In the 
notation \tup{\ref{not1.1}} with $\vec{Q} = \vec{\Gamma}$ these 
triangles are
\[ M_{\tau(y)} \xrightarrow{(g_{i, j})}
\bigoplus_{i = 1}^{m} \bigoplus_{j = 1}^{d_{i}} 
M_{x_{i}} \xrightarrow{(f_{i, j})^{\mrm{t}}} 
M_{y} \to M_{\tau(y)}[1] . \]
\item A morphism
$(g_{i, j}) : M_{\tau(y)} \to 
\bigoplus_{i = 1}^{m} \bigoplus_{j = 1}^{d_{i}} M_{x_{i}}$
is a source morphism iff 
for all $i$, $\{ g_{i, j} \}_{j = 1}^{d_{i}}$ is a basis of 
$\opn{Irr}_{\cat{D}^{\mrm{b}}(\cat{mod} A)}(M_{\tau(y)},
M_{x_{i}})$. Likewise a morphism
${(f_{i, j})^{\mrm{t}}} : 
\bigoplus_{i = 1}^{m} \bigoplus_{j = 1}^{d_{i}} M_{x_{i}} 
\to M_{y}$ is a sink morphism iff for all $i$,
$\{ f_{i, j} \}_{j = 1}^{d_{i}}$ is  basis of 
$\opn{Irr}_{\cat{D}^{\mrm{b}}(\cat{mod} A)}(M_{x_{i}},
M_{y})$. 
\end{enumerate}
\end{thm}

\begin{proof}
1. This is implicit in \cite{Ha} Sections I.4 and I.5. In 
particular \cite{Ha} Lemma I.5.2 shows that for any indecomposable 
object $M \in \cat{D}^{\mrm{b}}(\cat{mod} A)$ the ring \newline 
$\opn{End}_{\cat{D}^{\mrm{b}}(\cat{mod} A)}(M)$ is local.

\medskip \noindent
2. See \cite{Ha} Corollary I.4.9.

\medskip \noindent
3. According to \cite{Ha} Theorem I.4.6 and Lemma I.4.8, for each 
$y \in \vec{\Gamma}_{0}$ there exists such an Auslander-Reiten 
triangle. By \cite{Ha} Proposition I.4.3 these are all the 
Auslander-Reiten triangles, up to isomorphism. 

\medskip \noindent
4. Since source and sink morphism depend only on the structure of 
$k$-linear additive category on $\cat{D}^{\mrm{b}}(\cat{mod} A)$
(cf.\ \cite{Ha} Section I.4.5) we may use \cite{Rl} Lemma 2.2.3.
\end{proof}

The Auslander-Reiten quiver 
$\vec{\Gamma}(\cat{D}^{\mrm{b}}(\cat{mod} A))$ 
contains the quiver $\vec{\Delta}$, as the
full subquiver with vertices corresponding to the indecomposable
projective $A$-modules, under the inclusion
$\cat{mod} A \subset \cat{D}^{\mrm{b}}(\cat{mod} A)$.

\begin{dfn} \label{dfn2.1}
We call a connected component of 
$\vec{\Gamma}(\cat{D}^{\mrm{b}}(\cat{mod} A))$
{\em irregular} if it is isomorphic to the connected component 
containing $\vec{\Delta}$, and we denote by
$\vec{\Gamma}^{\mrm{irr}}$
the disjoint union of all irregular components of 
$\vec{\Gamma}(\cat{D}^{\mrm{b}}(\cat{mod} A))$.
\end{dfn}

The name ``irregular'' is inspired by \cite{ARS} Section VIII.4,
where regular components of $\vec{\Gamma}(\cat{mod} A)$ are 
discussed. The quiver $\vec{\Gamma}^{\mrm{irr}}$ will be of special 
interest to us. It's structure is explained in Theorem \ref{thm2.2} 
below. But first we need to recall the following definition due to
Riedtmann \cite{Rn},

\begin{dfn} \label{dfn2.2}
 From the quiver $\vec{\Delta}$ one can construct another quiver,
denoted by $\vec{\mbb{Z}} \vec{\Delta}$. The vertex set of
$\vec{\mbb{Z}} \vec{\Delta}$ is $\mbb{Z} \times \vec{\Delta}_{0}$, 
and for every arrow $x \xrightarrow{\alpha} y$ in $\vec{\Delta}$ 
there are arrows
$(n, x) \xrightarrow{(n, \alpha)} (n, y)$
and
$(n, y) \xrightarrow{(n, \alpha^{*})} (n + 1, x)$
in $\vec{\mbb{Z}} \vec{\Delta}$. 
\end{dfn}

The function $\tau (n, x) = (n - 1, x)$
makes $\vec{\mbb{Z}} \vec{\Delta}$ into a stable translation quiver.
Observe that $\tau$ is an automorphism of
the quiver $\vec{\mbb{Z}} \vec{\Delta}$, not just of the vertex set
$(\vec{\mbb{Z}} \vec{\Delta})_{0}$.
$\vec{\mbb{Z}} \vec{\Delta}$ is equipped with a polarization $\mu$,
given by
$\mu(n + 1, \alpha) = (n, \alpha^{*})$
and
$\mu(n, \alpha^{*}) = (n, \alpha)$.
See Figures \ref{fig2} and \ref{fig3} in Section 5 for examples.
We identify $\vec{\Delta}$ with the subquiver
$\{ 0 \} \times \vec{\Delta} \subset \vec{\mbb{Z}}\vec{\Delta}$.

Next let us define a quiver
$\mbb{Z} \times (\vec{\mbb{Z}} \vec{\Delta}) :=
\coprod_{m \in \mbb{Z}} \vec{\mbb{Z}} \vec{\Delta}$;
the connected components are
$\{ m \} \times (\vec{\mbb{Z}} \vec{\Delta})$, $m \in \mbb{Z}$.
Define an automorphism $\sigma$ of
$\mbb{Z} \times (\vec{\mbb{Z}} \vec{\Delta})$
by the action $\sigma(m) = m + 1$ on the first factor.
There is a translation $\tau$ and a polarization $\mu$ 
of $\mbb{Z} \times (\vec{\mbb{Z}} \vec{\Delta})$
that extend those of
$\vec{\mbb{Z}} \vec{\Delta} \cong \{ 0 \}  \times 
(\vec{\mbb{Z}} \vec{\Delta})$
and commute with $\sigma$.

The auto-equivalences $\sigma$ and $\tau$ of 
$\cat{D}^{\mrm{b}}(\cat{mod} A)$ induce commuting
permutations of $\vec{\Gamma}_{0}$, which we also denote by $\sigma$ 
and $\tau$ respectively. 

\begin{thm}[Happel] \label{thm2.2}
\begin{enumerate}
\item If $A$ has finite representation type then there is a
unique isomorphism of quivers
\[ \rho : \vec{\mbb{Z}} \vec{\Delta} \iso
\vec{\Gamma}^{\mrm{irr}} \]
which is the identity on $\vec{\Delta}$ and commutes with $\tau$
on vertices. Furthermore 
$\vec{\Gamma}^{\mrm{irr}} = 
\vec{\Gamma}(\cat{D}^{\mrm{b}}(\cat{mod} A))$.
\item If $A$ has infinite representation type
then there exists an isomorphism of quivers
\[ \rho : \mbb{Z} \times (\vec{\mbb{Z}} \vec{\Delta}) \iso
\vec{\Gamma}^{\mrm{irr}} \]
which is the identity on $\vec{\Delta}$ and commutes with $\tau$
and $\sigma$ on vertices.
If $\vec{\Delta}$ is a tree then the isomorphism $\rho$
is unique.
\end{enumerate}
\end{thm}

\begin{proof}
This is essentially \cite{Ha} Proposition I.5.5 and Corollary I.5.6.
\end{proof}

Fix once and for all for every vertex 
$x \in \vec{\Gamma}^{\mrm{irr}}_{0}$
an indecomposable object
$M_{x} \in \cat{D}^{\mrm{b}}(\cat{mod} A)$
which represents $x$, and such that $M_{x} = P_{x}$ for
$x \in \vec{\Delta}_{0}$.
Define $\msf{B} \subset \cat{D}^{\mrm{b}}(\cat{mod} A)$
to be the full subcategory with objects
$\{ M_{x} \mid x \in (\vec{\mbb{Z}} \vec{\Delta})_{0} \}$.

The additive category $\msf{B} \times \mbb{N}$ is also 
Krull-Schmidt, so for $M_{x}, M_{y} \in \msf{B}$ the two 
$k$-modules 
$\opn{Irr}_{\msf{B} \times \mbb{N}}(M_{x}, M_{y})$
and 
$\opn{Irr}_{\cat{D}^{\mrm{b}}(\cat{mod} A)}(M_{x}, M_{y})$
could conceivably differ (cf.\ \cite{Rl} Section 2.2). But this 
is not the case as we see in the lemma below. 

\begin{lem} \label{lem2.2}
Suppose $I \subset \mbb{Z}$ is a segment \tup{(}i.e.\ 
$I = \{i \in \mbb{Z} \mid a \leq i \leq b \}$ with 
$a, b \in \mbb{Z} \cup \{\pm \infty\}$\tup{)}. Let 
$\msf{B}(I) \subset \cat{D}^{\mrm{b}}(\cat{mod} A)$ 
be the full subcategory on the 
objects 
$M_{x}$, $x \in I \times \vec{\Delta}_{0} \subset 
\vec{\Gamma}(\cat{D}^{\mrm{b}}(\cat{mod} A))_{0}$. 
Then for any 
$M_{x}, M_{y} \in \msf{B}(I)$ one has
\[ \opn{Irr}_{\msf{B}(I) \times \mbb{N}}(M_{x}, M_{y}) \cong 
\opn{Irr}_{\cat{D}^{\mrm{b}}(\cat{mod} A)}(M_{x}, M_{y}) . \]
\end{lem}

\begin{proof}
Consider a sink morphism in $\cat{D}^{\mrm{b}}(\cat{mod} A)$ 
ending in $M_{(n, y)}$, 
$(n, y) \in (\vec{\mbb{Z}} \vec{\Delta})_{0}$. 
By Theorem \ref{thm2.1}(3) and Theorem \ref{thm2.2}, it is of the 
form 
${(f_{i, j})^{\mrm{t}}} : 
\bigoplus_{i = 1}^{m} \bigoplus_{j = 1}^{d_{i}} 
M_{(n - \epsilon_{i}, x_{i})}$ \linebreak
$\to M_{(n, y)}$ 
with $\epsilon_{i} \in \{ 0, 1 \}$ (cf.\ Notation \ref{not1.1}). 
 From the definition of a sink morphism we see that this is also a 
sink morphism in the category $\msf{B} \times \mbb{N}$. 

According to \cite{Rl} Lemma 2.2.3 (dual form), both $k$-modules 
\linebreak
$\opn{Irr}_{\msf{B} \times \mbb{N}}(M_{(n - \epsilon_{i}, x_{i})}, 
M_{(n, y)})$
and
$\opn{Irr}_{\cat{D}^{\mrm{b}}(\cat{mod} A)}
(M_{(n - \epsilon_{i}, x_{i})}, M_{(n, y)})$
have the morphisms $f_{i, 1}, \ldots, f_{i, d_{i}}$ 
as basis. And there are no irreducible morphisms $N \to M_{(n, y)}$ 
for indecomposable objects $N$ not isomorphic to one of the
$M_{(n - \epsilon_{i}, x_{i})}$, in either category. Thus the lemma is 
proved for $\msf{B}(I) = \msf{B}$.

Let $x, y \in \vec{\Delta}_{0}$ and $l, n \in \mbb{Z}$. If 
$\opn{Hom}(M_{(l, x)}, M_{(n, y)}) \neq 0$ then necessarily 
$l \leq n$. This is clear for $l = 0$, since  $M_{(0, x)}$ is a 
projective module, and an easy calculation shows that for $n < 0$,
\[ \mrm{H}^{0}(M_{(n, y)}) \cong
\mrm{H}^{0}(A^{*}[-1] \otimes^{\mrm{L}}_{A} \cdots 
\otimes^{\mrm{L}}_{A} A^{*}[-1] \otimes^{\mrm{L}}_{A} 
M_{(0, y)}) = 0 . \]
In general we can translate by $\tau^{-l}$. 

Now take an arbitrary segment $I$. The paragraph above implies that 
for $n, l \in I$ and $i \geq 0$, 
$\opn{rad}^{i}_{\msf{B}(I) \times \mbb{N}}(M_{(l, x)}, M_{(n, y)}) =
\opn{rad}^{i}_{\msf{B} \times \mbb{N}}(M_{(l, x)}, M_{(n, y)})$.
Hence \linebreak
$\opn{Irr}_{\msf{B}(I) \times \mbb{N}}(M_{(l, x)}, M_{(n, y)}) =
\opn{Irr}_{\msf{B} \times \mbb{N}}(M_{(l, x)}, M_{(n, y)})$.
\end{proof}

Henceforth we shall simply write $\opn{Irr}(M_{x}, M_{y})$ when 
$x, y \in (\vec{\mbb{Z}} \vec{\Delta})_{0}$. 
The lemma implies that the quiver of the category $\msf{B}(I)$ is 
the full subquiver 
$\vec{I} \vec{\Delta} \subset \vec{\mbb{Z}} \vec{\Delta}$.

Note that for $I = \{ 0 \}$ we get $\msf{B}(I) = \msf{P}$. 
Since $\msf{P}$ is canonically equivalent to 
$k \bra{\vec{\Delta}}$, 
there is a full faithful $k$-linear functor
$G_{0} : k \bra{\vec{\Delta}} \to \msf{B}$
such that $G_{0} x = M_{x} = P_{x}$ for every vertex
$x \in \vec{\Delta}_{0}$,
and 
$\{ G_{0}(\alpha_{j}) \}_{j = 1}^{d(x, y)}$
is a basis of $\opn{Irr}(M_{x}, M_{y})$
for every pair of vertices $x, y$, where
$\alpha_{1}, \cdots, \alpha_{d(x, y)}$ are the arrows 
$\alpha_{j} : x \to y$.

\begin{thm} \label{thm2.3}
Let $\vec{\Delta}$ be a finite quiver without oriented cycles, 
$A = k \vec{\Delta}$ its path algebra, 
$k \bra{\vec{\mbb{Z}} \vec{\Delta}, I_{\mrm{m}}}$
the mesh category \tup{(}Definitions \tup{\ref{dfn2.2}} and 
\tup{\ref{dfn1.1})} and 
$\msf{B} \subset \cat{D}^{\mrm{b}}(\cat{mod} A)$
the full subcategory on the objects
$\{ M_{x} \}_{x \in (\vec{\mbb{Z}} \vec{\Delta})_{0}}$. 
Then there is a $k$-linear functor
\[ G : k \bra{\vec{\mbb{Z}} \vec{\Delta}, I_{\mrm{m}}} 
\to \msf{B} \]
such that
\begin{enumerate}
\rmitem{i} $G x = M_{x}$ for each vertex
$x \in (\vec{\mbb{Z}} \vec{\Delta})_{0}$.
\rmitem{ii} $G|_{k \bra{\vec{\Delta}}} = G_{0}$.
\rmitem{iii} $G$ is full and faithful.
\end{enumerate}
Moreover, the functor $G$ is unique up to isomorphism.
\end{thm}

In other words, there is a unique equivalence $G$ extending $G_{0}$.

\begin{proof}
Let
$\vec{Q}^{+} \subset  \vec{\mbb{Z}} \vec{\Delta}$
be the full subquiver with vertex set
$\{ (n, y) \mid n \geq 0 \}$.
Given a vertex $(n, y)$ in $\vec{Q}^{+}$, denote by
$p(n,y)$ the number of its predecessors, i.e.\ the number of
vertices $(m, x)$ such that there is a path
$(m, x) \to \cdots \to (n,y)$ in $\vec{Q}^{+}$.
For any $p \geq 0$ let $\vec{Q}^{+}_{p}$ be the full subquiver
with vertex set
$\{ (n, y) \mid n \geq 0,\ p(n, y) \leq p \}$.
$\vec{Q}^{+}_{p}$ is a translation quiver with polarization, and 
$k \bra{\vec{Q}^{+}_{p}, I_{\mrm{m}}} \subset 
k \bra{\vec{\mbb{Z}} \vec{\Delta}, I_{\mrm{m}}}$
is a full subcategory.

By recursion on $p$, we will define a functor
$G : k \bra{\vec{Q}^{+}_{p}, I_{\mrm{m}}} \to \msf{B}$
satisfying conditions (i), (ii) and 
\begin{enumerate}
\rmitem{iv}  Let $x, y$ be a pair of vertices and let
$\alpha_{1}, \cdots, \alpha_{d(x, y)}$ be the arrows 
$\alpha_{j} : x \to y$. Then 
$\{ G(\alpha_{j}) \}_{j = 1}^{d(x, y)}$
is a basis of $\opn{Irr}(M_{x}, M_{y})$.
\end{enumerate}

Take $p \geq 0$.
It suffices to define $G(\alpha)$ for an arrow
$\alpha$ in $\vec{Q}^{+}_{p}$. These arrows fall into
three cases, according to their end vertex $(n, y)$:

\medskip \noindent
(a) $p(n, y) < p$, in which case any arrow $\alpha$ ending in
$(n, y)$ is in
$\vec{Q}^{+}_{p - 1}$, and $G(\alpha)$ is already
defined.

\medskip \noindent
(b) $p(n, y) = p$ and $n = 0$. Any arrow $\alpha$ ending in
$(n, y)$ is in $\vec{\Delta}$, so we define
$G(\alpha) := G_{0}(\alpha)$. By Lemma \ref{lem2.2} condition (iv) 
holds.

\medskip \noindent
(c) $p(n, y) = p$ and $n \geq 1$. In this case $(n, y)$ is a 
non-projective vertex in $\vec{Q}^{+}_{p}$, and we
consider the mesh ending at $(n, y)$. The vertices with arrows to
$(n, y)$ are $(n - \epsilon_{i}, x_{i})$,
where $i = 1, \ldots, m$; $x_{i} \in \vec{\Delta}_{0}$
and $\epsilon_{i} = 0, 1$
(cf.\ Notation \ref{not1.1}). 
Since 
$p(n - 1, y) < p(n - \epsilon_{i}, x_{i}) < p$
the arrows $\alpha_{i, j}$ are all in the quiver 
$\vec{Q}^{+}_{p - 1}$,
and hence $G(\alpha_{i, j})$ are defined.

According to condition (iv), Lemma  \ref{lem2.2} and Theorem 
\ref{thm2.1}(4) it follows that there exists an Auslander-Reiten 
triangle
\begin{equation} \label{eqn2.7}
M_{(n - 1, y)} \xrightarrow{(G(\alpha_{i, j}))}
\bigoplus_{i = 1}^{m} \bigoplus_{j = 1}^{d_{i}}
M_{(n - \epsilon_{i}, x_{i})}
\xrightarrow{(f_{i, j})^{\mrm{t}}}
M_{(n, y)} \to M_{(n - 1, y)}[1] 
\end{equation}
in $\cat{D}^{\mrm{b}}(\cat{mod} A)$. Define
\[ G(\beta_{i, j}) := f_{i, j} : M_{(n - \epsilon_{i}, x_{i})}
\to M_{(n, y)} . \]
Note that the mesh relation 
$\sum \beta_{i, j} \alpha_{i, j}$ in $k \bra{\vec{Q}^{+}_{p}}$ 
is sent by $G$ to
$\sum G(\beta_{i, j}) G(\alpha_{i, j}) = 0$, 
so we indeed have a functor
$G : k \bra{\vec{Q}^{+}_{p}, I_{\mrm{m}}} \to \msf{B}$.
Also, by Theorem \ref{thm2.1}(4), for any $i$ the set
$\{ G(\beta_{i, j}) \}_{j = 1}^{d_{i}}$ 
is a basis of 
$\opn{Irr}(M_{(n - \epsilon_{i}, x_{i})}, M_{(n, y)})$. 

Thus we obtain a functor 
$G : k \bra{\vec{Q}^{+}, I_{\mrm{m}}} \to \msf{B}$.

By symmetry we construct a functor
$G : k \bra{\vec{Q}^{-}, I_{\mrm{m}}} \to \msf{B}$
for negative vertices (i.e.\ $n \leq 0$), extending $G_{0}$. 
Putting the two together we obtain a functor
$G : k \bra{\vec{\mbb{Z}} \vec{\Delta}, I_{\mrm{m}}} \to \msf{B}$ 
satisfying conditions (i), (ii) and (iv). 

Let us prove $G$ is fully faithful.
For any $n \in \mbb{Z}$ there is a full subquiver
$\vec{\mbb{Z}}_{\geq n} \vec{\Delta} \subset 
\vec{\mbb{Z}} \vec{\Delta}$,
on the vertex set
$\{ (i, x) \mid i \geq n \}$.
Correspondingly there are full subcategories
$k \bra{\vec{\mbb{Z}}_{\geq n} \vec{\Delta}, I_{\mrm{m}}} \subset
k \bra{\vec{\mbb{Z}} \vec{\Delta}, I_{\mrm{m}}}$
and
$\msf{B}(\mbb{Z}_{\geq n}) \subset \msf{B}$.
It suffices to prove that
$G : k \bra{\vec{\mbb{Z}}_{\geq n} \vec{\Delta}, I_{\mrm{m}}} \to
\msf{B}(\mbb{Z}_{\geq n})$
is fully faithful. By Lemma \ref{lem2.2} the quiver of 
$\msf{B}(\mbb{Z}_{\geq n})$ is $\vec{\mbb{Z}}_{\geq n} \vec{\Delta}$, 
which is pre-projective. So we can use
the last two paragraphs in the proof of \cite{Rl} Lemma 2.3.3 almost 
verbatim.

Finally we shall prove that $G$ is unique up to isomorphism. 
Suppose 
$G' :  k \bra{\vec{\mbb{Z}} \vec{\Delta}, I_{\mrm{m}}} \to \msf{B}$
is another $k$-linear functor satisfying conditions (i)-(iii). 
We will show there is an isomorphism 
$\phi : G \iso G'$ that is the identity on $k \bra{\vec{\Delta}}$. 

By recursion on $p$ we shall exhibit an isomorphism 
$\phi : G|_{k \bra{\vec{Q}^{+}_{p}, I_{\mrm{m}}}}
\iso G'|_{k \bra{\vec{Q}^{+}_{p}, I_{\mrm{m}}}}$.
It suffices to consider case (c) above, so let $(n, y)$ be such a 
vertex. Then, because
$G'(\alpha_{i, j}) =  
\phi_{(n - \epsilon_{i}, x_{i})} G(\alpha_{i, j})
\phi_{(n - 1, y)}^{-1}$, 
we have
\[ \sum_{i, j} G'(\beta_{i, j}) \phi_{(n - \epsilon_{i}, x_{i})}
G(\alpha_{i, j}) = 
G' \bigl( \sum_{i, j} \beta_{i, j} \alpha_{i, j} \bigr) 
\phi_{(n - 1, y)} = 0 . \]
Applying $\opn{Hom}(-, M_{(n, y)})$ to the triangle (\ref{eqn2.7}) 
we obtain a morphism 
$a \in$ \linebreak
$\opn{End}(M_{(n, y)})$ 
such that 
$G'(\beta_{i, j}) \phi_{(n - \epsilon_{i}, x_{i})} = 
a G(\beta_{i, j})$.
Because $G'$ is faithful we see that $a \neq 0$, and since 
$\opn{End}(M_{(n, y)}) \cong k$ 
it follows that $a$ is invertible. Set
$\phi_{(n, y)} := a \in \opn{Aut}(M_{(n, y)})$.
This yields the desired isomorphism
$\phi : G|_{k \bra{\vec{Q}^{+}_{p}, I_{\mrm{m}}}}
\iso G'|_{k \bra{\vec{Q}^{+}_{p}, I_{\mrm{m}}}}$.

By symmetry the isomorphism $\phi$ extends to $\vec{Q}^{-}$.

\end{proof}

The uniqueness of $G$ gives the next corollary. 

\begin{cor} \label{cor2.9}
Let $F$ be a $k$-linear auto-equivalence of
$k \bra{\vec{\mbb{Z}} \vec{\Delta}, I_{\mrm{m}}}$
fixing all objects, and such that 
$F|_{k \bra{\vec{\Delta}}} \cong \mbf{1}_{k \bra{\vec{\Delta}}}$.
Then 
$F \cong \mbf{1}_{k \bra{\vec{\mbb{Z}} \vec{\Delta}, I_{\mrm{m}}}}$.
\end{cor}

\begin{rem} \label{rem2.2}
Beware that if $A$ has infinite representation type then 
$k \bra{\vec{\Gamma}^{\mrm{irr}}, I_{\mrm{m}}}$ is {\em not} 
equivalent to the full subcategory of 
$\cat{D}^{\mrm{b}}(\cat{mod} A)$ on the objects
$\{ M_{x} \}_{x \in \vec{\Gamma}^{\mrm{irr}}_{0}}$. 
This is because there are nonzero morphisms from the projective 
modules (vertices in the component 
$\vec{\mbb{Z}} \vec{\Delta}$) to the injective modules 
(vertices in $\{ 1 \} \times \vec{\mbb{Z}} \vec{\Delta}$).
\end{rem}

\section{The Representation of $\opn{DPic}_{k}(A)$ on the Quiver
$\vec{\Gamma}^{\mrm{irr}}$}

This section contains the proof of the main result of the paper, 
Theorem \ref{thm0.2} (restated here as Theorem \ref{thm3.2}). It is 
deduced from the more technical Theorem \ref{thm3.1}. 
Throughout $k$ is an algebraically closed field, $\vec{\Delta}$ is a 
connected finite quiver without oriented cycles, and 
$A = k \vec{\Delta}$ is the path algebra. We use the notation of 
previous sections.

Recall that $A^{*} = \opn{Hom}_{k}(A, k)$ is a tilting complex. 
We shall denote by $\tau$ the class of $A^{*}[-1]$
in $\opn{DPic}_{k}(A)$, and by $\sigma$ the class of $A[1]$.
We identify an element $T \in \opn{DPic}_{k}(A)$ and the induced 
auto-equivalence $F = T \otimes^{\mrm{L}}_{A} -$ of 
$\cat{D}^{\mrm{b}}(\cat{mod} A)$.

\begin{lem} \label{lem3.8}
$\tau$  and $\sigma$ are in the center of
$\opn{DPic}_{k}(A)$.
\end{lem}

\begin{proof}
The fact that $\sigma$ is in the center of $\opn{DPic}_{k}(A)$
is trivial. As for $\tau$, this follows immediately from \cite{Rd} 
Proposition 5.2 (or by \cite{BO} Proposition 1.3, since 
$A^{*} \otimes^{\mrm{L}}_{A} -$ is 
the Serre functor of $\cat{D}^{\mrm{b}}(\cat{mod} A)$).
\end{proof}

In Definition \ref{dfn2.1} we introduced the quiver 
$\vec{\Gamma}^{\mrm{irr}}$. Recall that for a vertex 
$x \in \vec{\Gamma}^{\mrm{irr}}_{0}$, 
$M_{x} \in \cat{D}^{\mrm{b}}(\cat{mod} A)$ is the representative 
indecomposable object.

\begin{lem} \label{lem3.6}
There is a group homomorphism
\[ q : \opn{DPic}_{k}(A) \to
\opn{Aut}(\vec{\Gamma}^{\mrm{irr}}_{0}; d)^{\bra{\tau, \sigma}} \]
such that $q(F)(x) = y$ iff $F M_{x} \cong M_{y}$.
\end{lem}

\begin{proof}
Given an auto-equivalence $F$ of $\cat{D}^{\mrm{b}}(\cat{mod} A)$,
the formula $q(F)(x) = y$ iff $F M_{x} \cong M_{y}$ defines a 
permutation $q(F)$ of 
$\vec{\Gamma}_{0}(\cat{D}^{\mrm{b}}(\cat{mod} A))$
that preserves arrow-multiplicities. Hence it restricts to a 
permutation of $\vec{\Gamma}^{\mrm{irr}}_{0}$. 
By Lemma \ref{lem3.8}, $q(F)$ commutes with $\tau$ and $\sigma$. 
\end{proof}

The group $\opn{Out}^{0}_{k}(A)$ was defined to be the identity  
component of $\opn{Out}_{k}(A)$.

\begin{lem} \label{lem3.2}
$\opn{Ker}(q) = \opn{Out}^{0}_{k}(A)$. 
\end{lem}

\begin{proof}
Let $T \in \opn{DPic}_{k}(A)$. By Theorem \ref{thm2.2} we know
that 
$\vec{\Gamma}^{\mrm{irr}}_{0} = 
\bigcup_{i, j \in \mbb{Z}} \tau^{i} \sigma^{j}(\vec{\Delta}_{0})$.
Hence by Lemma \ref{lem3.8}, $T \in \opn{Ker}(q)$ iff $T$ acts 
trivially on the set $\vec{\Delta}_{0}$. In particular we see that 
$\opn{Ker}(q) \subset \opn{Pic}_{k}(A)$. Now use Proposition 
\ref{prop1.1}.
\end{proof}

\begin{lem} \label{lem3.11}
Suppose $A$ has finite representation type. Then
$\sigma$ is in the center of the group
$\opn{Aut}(\vec{\Gamma}^{\mrm{irr}})^{\bra{\tau}}$.
\end{lem}

\begin{proof}
According to \cite{Rn} Section 2, the 
group $\opn{Aut}(\vec{\mbb{Z}} \vec{\Delta})^{\bra{\tau}}$
is abelian in all cases except $D_{4}$. But a direct calculation
in this case  (cf.\ Theorem \ref{thm4.1}) 
gives $\sigma = \tau^{-3}$.
\end{proof}

Before we can talk about the mesh category 
$k \bra{\vec{\Gamma}^{\mrm{irr}}, I_{\mrm{m}}}$
of the quiver
$\vec{\Gamma}^{\mrm{irr}}$, we have to fix a polarization $\mu$ on 
it. If the quiver $\vec{\Delta}$ has no multiple arrows then so 
does $\vec{\Gamma}^{\mrm{irr}}$ (by Theorem \ref{thm2.2}), and hence
there is a unique polarization on it. If $\vec{\Delta}$ isn't a tree 
let us choose an isomorphism 
$\rho : \mbb{Z} \times (\vec{\mbb{Z}} \vec{\Delta}) \iso
\vec{\Gamma}^{\mrm{irr}}$
as in that theorem. This determines a polarization $\mu$ on 
$\vec{\Gamma}^{\mrm{irr}}$.
We also get a lifting of the permutation $\sigma$ to an
auto-equivalence of 
$k \bra{\vec{\Gamma}^{\mrm{irr}}, I_{\mrm{m}}}$.

\begin{lem} \label{lem3.12}
There are group homomorphisms
\[ p : \opn{Out}_{k}(k \bra{\vec{\Gamma}^{\mrm{irr}}, I_{\mrm{m}}})
\to \opn{Aut}(\vec{\Gamma}^{\mrm{irr}}_{0}; d)^{\bra{\tau}}  \]
and 
\[ r : \opn{Aut}(\vec{\Gamma}^{\mrm{irr}}_{0}; d)^{\bra{\tau}}
\to
\opn{Out}_{k}(k \bra{\vec{\Gamma}^{\mrm{irr}}, I_{\mrm{m}}}) \]
satisfying $p(F)(x) = F x$  for an auto-equivalence $F$ and a 
vertex $x$; $p r = 1$; and both $p$ and $r$ commute with $\sigma$.
\end{lem}

\begin{proof}
Since 
$\vec{\Gamma}^{\mrm{irr}} \cong 
\vec{\Gamma}(k \bra{\vec{\Gamma}^{\mrm{irr}}, I_{\mrm{m}}})$
we get a permutation 
$p(F) \in \opn{Aut}(\vec{\Gamma}^{\mrm{irr}}_{0}; d)$. 
Let's prove that $p(F)$ commutes with $\tau$ in 
$\opn{Aut}(\vec{\Gamma}^{\mrm{irr}}_{0})$. Consider a 
vertex $y \in \vec{\Gamma}^{\mrm{irr}}_{0}$. In the Notation 
\ref{not1.1}, there are vertices $x_{i}$ and irreducible morphisms
$\{ F(\alpha_{i, j}) \}_{j=1}^{d_{i}}$ and
$\{ F(\beta_{i, j}) \}_{j=1}^{d_{i}}$
that form bases of 
$\opn{Irr}_{k \bra{\vec{\Gamma}^{\mrm{irr}}, I_{\mrm{m}}}}
(F \tau y, F x_{i})$
and
$\opn{Irr}_{k \bra{\vec{\Gamma}^{\mrm{irr}}, I_{\mrm{m}}}}
(F x_{i}, F y)$
respectively. Since we have 
\[ \sum F(\beta_{i, j}) F(\alpha_{i, j}) = 0 \in
\opn{rad}_{k \bra{\vec{\Gamma}^{\mrm{irr}}, I_{\mrm{m}}}}^{2}
(F \tau y, F y) / 
\opn{rad}_{k \bra{\vec{\Gamma}^{\mrm{irr}}, I_{\mrm{m}}}}^{3}
(F \tau y, F y) \]
this must be a multiple of a mesh relation. Hence 
$F \tau y =  \tau F y$.

Finally to define $r$ we have to split 
$\opn{Aut}(\vec{\Gamma}^{\mrm{irr}}) \surj
\opn{Aut}(\vec{\Gamma}^{\mrm{irr}}_{0}; d)$
consistently with $\mu$. It suffices to order the set of 
arrows $\{ \alpha: x \to y \}$ for every pair of vertices 
$x, y \in \vec{\Gamma}^{\mrm{irr}}_{0}$
consistently with $\mu$. We only have to worry about this when 
$A$ has infinite representation type. For any 
$x, y \in \vec{\Delta}_{0}$ choose some ordering of the set
$\{ \alpha: x \to y \}$. Using $\mu$ and $\sigma$ this ordering can 
be transported to all of 
$\mbb{Z} \times (\vec{\mbb{Z}} \vec{\Delta})$.
By the isomorphism $\rho$ of Theorem \ref{thm2.2} the ordering is 
copied to $\vec{\Gamma}^{\mrm{irr}}$.
\end{proof}

\begin{lem} \label{lem3.5}
There exists a group homomorphism
\[ \tilde{q} : \opn{DPic}_{k}(A) \to
\opn{Out}_{k}(k \bra{\vec{\Gamma}^{\mrm{irr}}, I_{\mrm{m}}})  \]
such that
$p \tilde{q} = q$.
\end{lem}

\begin{proof}
Choose an equivalence
$G : k \bra{\vec{\mbb{Z}} \vec{\Delta}, I_{\mrm{m}}} \to \msf{B}$
as in Theorem \ref{thm2.3}. If $A$ has infinite representation 
type then the isomorphism $\rho$ we have chosen (as in Theorem 
\ref{thm2.2}) tells us how to extend $G$ to an equivalence
$G : k \bra{\vec{\Gamma}^{\mrm{irr}}, I_{\mrm{m}}} \to
\coprod_{l \in \mbb{Z}} \msf{B}[l]$
that commutes with $\sigma$ (cf.\ Remark \ref{rem2.2}). 

Let $F$ be a triangle auto-equivalence of  
$\cat{D}^{\mrm{b}}(\cat{mod} A)$.
Then $F$ induces a permutation $\pi = q(F) $ of the set
$\vec{\Gamma}^{\mrm{irr}}_{0}$
that commutes with $\sigma$. For every vertex
$x \in \vec{\Gamma}^{\mrm{irr}}_{0}$
choose an isomorphism
\[ \phi_{x} : F M_{x} \iso M_{\pi(x)} \]
in $\cat{D}^{\mrm{b}}(\cat{mod} A)$. Given an arrow
$\alpha : x \to y$ in $\vec{\Gamma}^{\mrm{irr}}$,
define the morphism
$\tilde{q}_{ \{ \phi_{x} \} }(F)(\alpha) : \pi(x) \to \pi(y)$
by the condition that the diagram
\[ \begin{CD}
F M_{x} @> F G(\alpha) >> F M_{y} \\
@V \phi_{x} VV  @V \phi_{y} VV  \\
M_{\pi(x)} @> G \tilde{q}_{ \{ \phi_{x} \} }(F)(\alpha) >> 
M_{\pi(y)}
\end{CD} \]
commutes. Then 
$\tilde{q}_{ \{ \phi_{x} \} }(F) \in 
\opn{Aut}_{k}(k \bra{\vec{\Gamma}^{\mrm{irr}}, I_{\mrm{m}}})$.

If $\{ \phi'_{x} \}$ is another choice of isomorphisms
$\phi'_{x} : F M_{x} \iso M_{\pi(x)}$
then $\{ \phi'_{x} \phi^{-1}_{x} \}$
is an isomorphism of functors
$\tilde{q}_{ \{ \phi_{x} \} }(F) \to
\tilde{q}_{ \{ \phi'_{x} \} }(F)$,
so the map
$\tilde{q} : \opn{DPic}_{k}(A) \to
\opn{Out}_{k}(k \bra{\vec{\Gamma}^{\mrm{irr}}, I_{\mrm{m}}})$
is independent of these choices.

It is easy to check that $\tilde{q}$ respects composition of
equivalences.
\end{proof}

\begin{thm} \label{thm3.1}
Let $A$ be an indecomposable basic hereditary finite dimensional 
$k$-algebra with quiver $\vec{\Delta}$. Then the homomorphism
$\tilde{q}$ of Lemma \tup{\ref{lem3.5}} induces an isomorphism 
of groups 
\[ \opn{DPic}_{k}(A) \cong
\opn{Out}_{k}(k \bra{\vec{\Gamma}^{\mrm{irr}}, I_{\mrm{m}}})
^{\bra{\sigma}} \cong
\begin{cases}
\opn{Out}_{k} \bigl( k \bra{\vec{\mbb{Z}} \vec{\Delta}, I_{\mrm{m}}} 
\bigr) \quad 
\parbox[t]{3cm}{\textup{if } A \textup{ has finite representation
type}} \\[6mm]
\opn{Out}_{k} \bigl( k \bra{\vec{\mbb{Z}} \vec{\Delta}, I_{\mrm{m}}}
\bigr) \times \bra{\sigma} \quad
\parbox[t]{2cm}{\textup{otherwise.}}
\end{cases} \]
\end{thm}

\begin{proof}
The proof has three parts.

\medskip \noindent 1. We show that the homomorphism
\[ \tilde{q} : \opn{DPic}_{k}(A) \to
\opn{Out}_{k}(k \bra{\vec{\Gamma}^{\mrm{irr}}, I_{\mrm{m}}})  \]
of Lemma \ref{lem3.5} is injective.
Let $T$ be a two-sided tilting complex such that
$\tilde{q}(T) \cong
\bsym{1}_{k \bra{\vec{\Gamma}^{\mrm{irr}}, I_{\mrm{m}}}}$. 
Then the permutation $q(T)$ fixes the vertices of
$\vec{\Delta} \subset \vec{\Gamma}^{\mrm{irr}}$.
Using the fact that
$A \cong \bigoplus_{x \in \vec{\Delta}_{0}} M_{x}$
we see that
$T \cong A$ in $\cat{D}^{\mrm{b}}(\cat{mod} A)$.
Replacing $T$ with $\opn{H}^{0} T$ we may assume
$T$ is a single bimodule.
According to \cite{Ye} Proposition 2.2, we see that $T$ is
actually an invertible bimodule. Since 
$k \bra{\vec{\Delta}} \to 
k \bra{\vec{\mbb{Z}} \vec{\Delta}, I_{\mrm{m}}}$
is full we get 
$\tilde{q}(T)|_{k \bra{\vec{\Delta}}} \cong
\bsym{1}_{k \bra{\vec{\Delta}}}$.
Hence by Morita theory we have $T \cong A$ as bimodules.

\medskip \noindent 2. Assume $A$ has finite representation type,
so that 
$\vec{\Gamma}^{\mrm{irr}} \cong \vec{\mbb{Z}} \vec{\Delta}$. 
We prove that
\[ \tilde{q} : \opn{DPic}_{k}(A) \to
\opn{Out}_{k}(k \bra{\vec{\mbb{Z}} \vec{\Delta}, I_{\mrm{m}}})  \]
is surjective.

Consider a $k$-linear auto-equivalence $F$ of
$k \bra{\vec{\mbb{Z}} \vec{\Delta}, I_{\mrm{m}}}$. Let
$\pi := p(F) \in$ \linebreak 
$\opn{Aut}((\vec{\mbb{Z}} \vec{\Delta})_{0}; d)^{\bra{\tau}}
\cong \opn{Aut}(\vec{\mbb{Z}} \vec{\Delta})^{\bra{\tau}}$
as in the proof of Lemma \ref{lem3.5}. According to Lemma 
\ref{lem3.11}, $\pi$ commutes with $\sigma$. Define
\[ M := \bigoplus_{x \in \vec{\Delta}_{0}} M_{\pi(0, x)}
\in \cat{D}^{\mrm{b}}(\cat{mod} A) . \]
Then for any $x, y \in \vec{\Delta}_{0}$ and integers $n, i$
the equivalence
$G : k \bra{\vec{\mbb{Z}} \vec{\Delta}, I_{\mrm{m}}} \to \msf{B}$
of Theorem \ref{thm2.3} produces isomorphisms
\[ \begin{aligned}
\opn{Hom}_{\cat{D}^{\mrm{b}}(\cat{mod} A)}
(M_{\pi(0, x)}, M_{(n, y)}[i])
& \cong
\opn{Hom}_{k \bra{\vec{\mbb{Z}} \vec{\Delta}, I_{\mrm{m}}}}
(\pi(0, x), \sigma^{i}(n, y)) \\
& \cong
\opn{Hom}_{k \bra{\vec{\mbb{Z}} \vec{\Delta}, I_{\mrm{m}}}}
((0, x), \sigma^{i} \pi^{-1}(n, y)) \\
& \cong
\opn{Hom}_{\cat{D}^{\mrm{b}}(\cat{mod} A)}
(M_{(0, x)}, M_{\pi^{-1}(n, y)}[i]) .
\end{aligned} \]

Therefore
\[ \begin{aligned}
\opn{Hom}_{\cat{D}^{\mrm{b}}(\cat{mod} A)}(M, M[i])
& \cong
\bigoplus_{x, y \in \vec{\Delta}_{0}}
\opn{Hom}_{\cat{D}^{\mrm{b}}
(\cat{mod} A)}(M_{(0, x)}, M_{(0, y)}[i]) \\
& \cong
\begin{cases}
A^{\circ} & \text{if } i = 0 \\
0 & \text{otherwise} .
\end{cases}
\end{aligned} \]
Also for any $(n, y)$ there is some integer $i$ and
$x \in \vec{\Delta}_{0}$ such that
\[ \opn{Hom}_{\cat{D}^{\mrm{b}}(\cat{mod} A)}
(M_{\pi(0, x)}, M_{(n, y)}[i])
\neq 0 . \]
Since any object $N \in \cat{D}^{\mrm{b}}(\cat{mod} A)$ is a direct
sum of indecomposables $M_{(n, y)}$, this implies that 
$\opn{R} \opn{Hom}_{A}(M, N) \neq 0$ if $N \neq 0$. By \cite{Ye} 
Theorem 1.8 and the proof of ``(ii) $\Rightarrow$ (i)'' 
of \cite{Ye} Theorem 1.6 there exists a two-sided tilting complex 
$T$ with $T \cong M$ in $\cat{D}(\cat{Mod} A)$ 
(cf.\ \cite{Rd} Section 3). Replacing 
$F$ with $\tilde{q}(T^{\vee}) F$, 
where $T^{\vee} := \opn{R} \opn{Hom}_{A}(T, A)$,
we can assume that $p(F)$ is trivial.

Now that $p(F)$ is trivial, $F$ restricts to an
auto-equivalence of $k \bra{\vec{\Delta}}$, and by Proposition 
\ref{prop1.1} we have 
$F|_{k \bra{\vec{\Delta}}} \cong \mbf{1} _{k \bra{\vec{\Delta}}}$. 
Then Corollary \ref{cor2.9} tells us 
$F \cong \mbf{1}_{k \bra{\vec{\mbb{Z}} \vec{\Delta}, I_{\mrm{m}}}}$.

\medskip \noindent 3. Assume $A$ has infinite representation type.
Then the quiver isomorphism $\rho$ of Theorem \ref{thm2.2} induces
a group isomorphism
\[ \opn{Out}_{k}(k \bra{\vec{\Gamma}^{\mrm{irr}}, I_{\mrm{m}}})
^{\bra{\sigma}}
\cong
\opn{Out}_{k}(k \bra{\vec{\mbb{Z}} \vec{\Delta}, I_{\mrm{m}}})
\times \bra{\sigma} , \]
and $\bra{\sigma} \cong \mbb{Z}$.
We prove that
\[ \tilde{q} : \opn{DPic}_{k}(A) \to
\opn{Out}_{k}(k \bra{\vec{\mbb{Z}} \vec{\Delta}, I_{\mrm{m}}})
\times \mbb{Z} \]
is surjective.

Take an auto-equivalence $F$ of
$k \bra{\vec{\mbb{Z}} \vec{\Delta}, I_{\mrm{m}}}$, and write
$\pi := p(F) \in$ \newline 
$\opn{Aut}((\vec{\mbb{Z}} \vec{\Delta})_{0}; d)^{\bra{\tau}}$.
After replacing $F$ with $\tau^{j} F$ for suitable
$j \in \mbb{Z}$, we can assume that
$\pi(0, x) \in \vec{\mbb{Z}}_{\geq 0} \vec{\Delta}$
for all $x \in \vec{\Delta}_{0}$.
Because $\vec{\mbb{Z}}_{\geq 0} \vec{\Delta}$ is the preprojective
component of $\vec{\Gamma}(\cat{mod} A)$ (cf.\ \cite{Rl}), we get 
\[ M := \bigoplus_{x \in \vec{\Delta}_{0}} M_{\pi(0, x)}
\in \cat{mod} A . \]
As in part 2 above, $\opn{End}_{A}(M) = A^{\circ}$.
Since $M$ is a complete slice, \cite{HR} Theorem 7.2 says that $M$ 
is a tilting module. So $M$ is a two-sided tilting complex over $A$.
Replacing $F$ by $\tilde{q}(M^{\vee}) F$ we can 
assume $p(F)$ is trivial. Let $P$ be an invertible bimodule such that 
$\tilde{q}(P)|_{k \bra{\vec{\Delta}}} \cong 
F|_{k \bra{\vec{\Delta}}}$.
Replacing $F$ with $\tilde{q}(P^{\vee}) F$
we get 
$F|_{k \bra{\vec{\Delta}}} \cong \mbf{1}_{k \bra{\vec{\Delta}}}$.
Then by Corollary \ref{cor2.9} we get 
$F \cong \mbf{1}_{\bra{\vec{\mbb{Z}} \vec{\Delta}, I_{\mrm{m}}}}$.
\end{proof}

The next theorem is Theorem \ref{thm0.2} in the Introduction.

\begin{thm} \label{thm3.2}
Let $A$ be an indecomposable basic hereditary finite dimensional 
algebra over an algebraically closed field $k$, with quiver 
$\vec{\Delta}$.
\begin{enumerate}
\item There is an exact sequence of groups
\[ 1 \to \opn{Out}^{0}_{k}(A) \to \opn{DPic}_{k}(A) 
\xrightarrow{q} 
\opn{Aut}(\vec{\Gamma}^{\mrm{irr}}_{0} ; d)^{\bra{\tau, \sigma}}
\to 1 . \]
This sequence splits.
\item If $A$ has finite representation type then there is an 
isomorphism of groups
\[ \opn{DPic}_{k}(A) \cong
\opn{Aut}(\vec{\mbb{Z}} \vec{\Delta})^{\bra{\tau}} . \]
\item If $A$ has infinite representation type then 
there is an isomorphism of groups
\[ \opn{DPic}_{k}(A) \cong
\bigl( \opn{Aut}((\vec{\mbb{Z}} \vec{\Delta})_{0}; d)^{\bra{\tau}}
\ltimes \opn{Out}^{0}_{k}(A) \bigr) \times \mbb{Z} . \]
\end{enumerate}
\end{thm}

\begin{proof}
1. By Theorem \ref{thm3.1} and Lemma \ref{lem3.12} the homomorphism $q$ 
is surjective. Lemma \ref{lem3.2} identifies $\opn{Ker}(q)$.

\medskip \noindent
2. If $A$ has finite representation type then $\vec{\Delta}$ is a 
tree, so $\opn{Out}_{k}^{0}(A) = 1$ by Proposition  \ref{prop1.1}.
By Theorem \ref{thm2.2} and Lemma \ref{lem3.11} we get
\[ \opn{Aut}(\vec{\Gamma}^{\mrm{irr}}_{0}; d)^{\bra{\tau, \sigma}}
\cong \opn{Aut}(\vec{\Gamma}^{\mrm{irr}})^{\bra{\tau, \sigma}}
\cong \opn{Aut}(\vec{\mbb{Z}} \vec{\Delta})^{\bra{\tau}} . \]

\medskip \noindent
3. If $A$ has infinite representation type then 
\[ \opn{Aut}(\vec{\Gamma}^{\mrm{irr}}_{0}; d)^{\bra{\tau, \sigma}} 
\cong 
\opn{Aut}((\vec{\mbb{Z}} \vec{\Delta})_{0}; d)^{\bra{\tau}} \times
\bra{\sigma} \]
by Theorem \ref{thm2.2}. We know that $\sigma$ is in the center 
of $\opn{DPic}_{k}(A)$.
\end{proof}

We end the section with the following problem.
 
\begin{prob} \label{prob3.1}
The Auslander-Reiten quiver 
$\vec{\Gamma}(\cat{D}^{\mrm{b}}(\cat{mod} A))$
is defined for any finite dimensional $k$-algebra $A$ of finite 
global dimension. Can the action of $\opn{DPic}_{k}(A)$ on 
$\vec{\Gamma}(\cat{D}^{\mrm{b}}(\cat{mod} A))$
be used to determine the structure of $\opn{DPic}_{k}(A)$ for any 
such $A$?
\end{prob} 

\section{Explicit Calculations}

In this section we calculate the group structure of 
$\opn{DPic}_{k}(A)$ for the path algebra $A = k \vec{\Delta}$
for several types of quivers. Throughout $S_{m}$ 
denotes the permutation group of $\{ 1, \ldots, m \}$. 

Suppose $\Delta$ is a tree. Given an orientation $\omega$ of the 
edge set $\Delta_{1}$, denote by $\vec{\Delta}_{\omega}$ the 
resulting quiver, and by $A_{\omega} := k \vec{\Delta}_{\omega}$.
If $\omega$ and $\omega'$ are two orientations of $\Delta$ then 
$\cat{D}^{\mrm{b}}(\cat{mod} A_{\omega}) \approx
\cat{D}^{\mrm{b}}(\cat{mod} A_{\omega'})$.
This equivalence will be discussed in the next section. For now we 
just note that the groups
$\opn{DPic}_{k}(A_{\omega}) \cong \opn{DPic}_{k}(A_{\omega'})$, 
so we are allowed to choose any orientation of $\Delta$
when computing these groups. This observation is relevant to 
Theorems \ref{thm4.1} and \ref{thm4.2} below.

\begin{figure}
\choosegraphics{
\[ \UseTips
\begin{array}{lr}
A_{n} \quad
\begin{xy} 
            (0,0)*+@{*}="1"*+!U{\scrp{1}},"1"
\ar@{->}    (8,0)*+@{*}="2"*+!U{\scrp{2}},"2" 
\ar@{->}    "2";(16,0)*+@{*}="3"*+!U{\scrp{3}},"3" 
\ar@{}     "3";(24,0)*+@{}="4","4" |*{\cdots}
\ar@{->}    "4";(32,0)*+@{*}="n"*+!U{\scrp{n}},"n" 
\end{xy} 
& \hspace{15mm} D_{n} \quad
\begin{xy} 
            (0,0)*+@{*}="1"*+!R{\scrp{1}},"1"
\ar@{->}    (0,8)*+@{*}="2"*+!R{\scrp{2}},"2" 
\ar@{->}    (0,-8)*+@{*}="3"*+!R{\scrp{3}},"3"
\ar@{->}    (8,0)*+@{*}="4"*+!U{\scrp{4}},"4" 
\ar@{->}    "4";(16,0)*+@{*}="5"*+!U{\scrp{5}},"5" 
\ar@{}      "5";(24,0)*+@{}="6","6" |*{\cdots}
\ar@{->}    "6";(32,0)*+@{*}="n"*+!U{\scrp{n}},"n" 
\end{xy} \\
E_{6} \quad
\begin{xy} 
            (0,0)*+@{*}="1"*+!U{\scrp{1}},"1"
\ar@{->}    (8,0)*+@{*}="2"*+!U{\scrp{2}},"2" 
\ar@{->}    "2";(16,0)*+@{*}="3"*+!U{\scrp{3}},"3" 
\ar@{<-}    "3";(24,0)*+@{*}="5"*+!U{\scrp{5}},"5" 
\ar@{<-}    "5";(32,0)*+@{*}="6"*+!U{\scrp{6}},"6" 
\ar@{<-}    "3";(16,8)*+@{*}="4"*+!R{\scrp{4}},"4" 
\end{xy} 
& \hspace{15mm} E_{7} \quad
\begin{xy} 
            (0,0)*+@{*}="1"*+!U{\scrp{1}},"1"
\ar@{->}    (8,0)*+@{*}="2"*+!U{\scrp{2}},"2" 
\ar@{->}    "2";(16,0)*+@{*}="3"*+!U{\scrp{3}},"3" 
\ar@{<-}    "3";(16,8)*+@{*}="4"*+!R{\scrp{4}},"4" 
\ar@{<-}    "3";(24,0)*+@{*}="5"*+!U{\scrp{5}},"5" 
\ar@{<-}    "5";(32,0)*+@{*}="6"*+!U{\scrp{6}},"6" 
\ar@{<-}    "6";(40,0)*+@{*}="7"*+!U{\scrp{7}},"7" 
\end{xy}  \\[5mm]
& E_{8} \quad
\begin{xy} 
            (0,0)*+@{*}="1"*+!U{\scrp{1}},"1"
\ar@{->}    (8,0)*+@{*}="2"*+!U{\scrp{2}},"2" 
\ar@{->}    "2";(16,0)*+@{*}="3"*+!U{\scrp{3}},"3" 
\ar@{<-}    "3";(16,8)*+@{*}="4"*+!R{\scrp{4}},"4" 
\ar@{<-}    "3";(24,0)*+@{*}="5"*+!U{\scrp{5}},"5" 
\ar@{<-}    "5";(32,0)*+@{*}="6"*+!U{\scrp{6}},"6" 
\ar@{<-}    "6";(40,0)*+@{*}="7"*+!U{\scrp{7}},"7" 
\ar@{<-}    "7";(48,0)*+@{*}="8"*+!U{\scrp{8}},"8" 
\end{xy}  
\end{array} \]
}{
\includegraphics[clip]{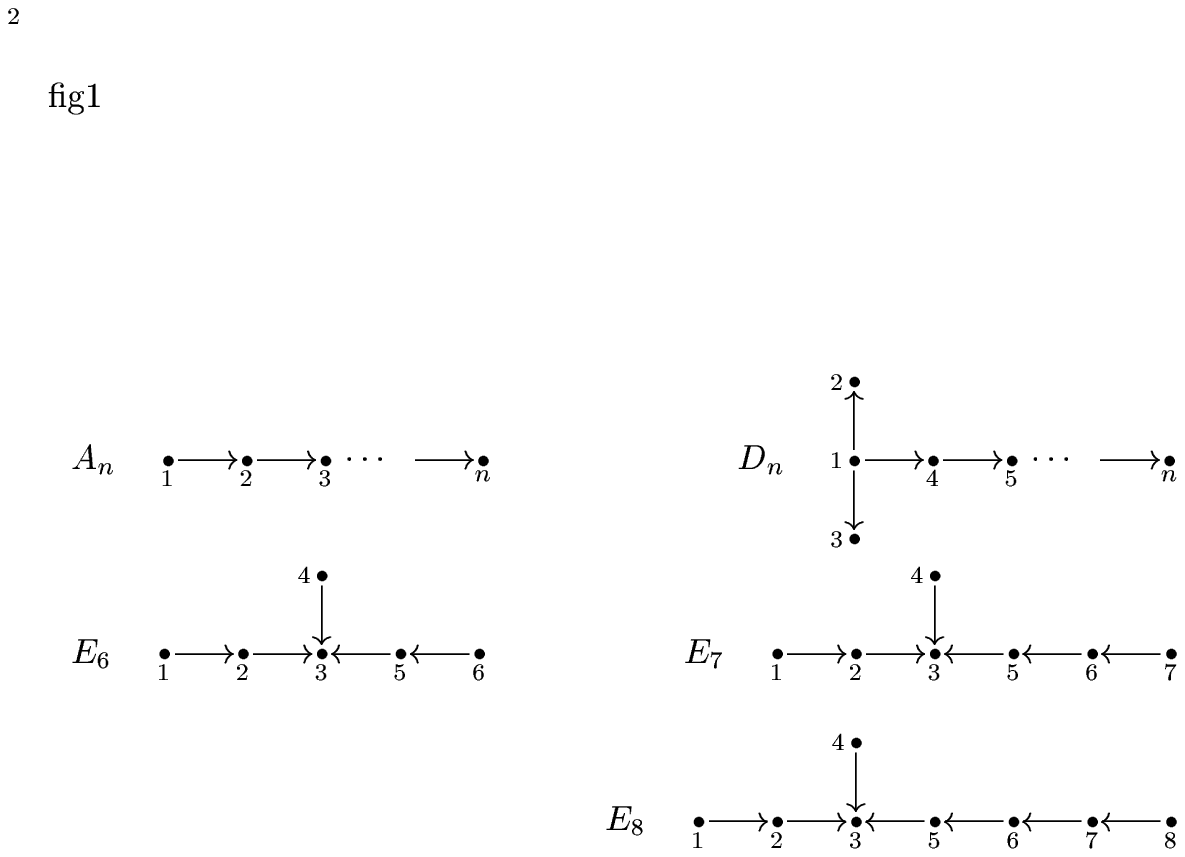}
}
\caption{Orientations for the Dynkin graphs} \label{fig1}
\end{figure}

\begin{thm} \label{thm4.1}
Let $\vec{\Delta}$ be a Dynkin quiver as shown in Figure 
\tup{\ref{fig1}}, and let $A := k \vec{\Delta}$ be the path algebra. 
Then 
$\opn{Pic}_{k}(A) \cong \opn{Aut}(\vec{\Delta})$
and
$\opn{DPic}_{k}(A) \cong
\opn{Aut}(\vec{\mbb{Z}} \vec{\Delta})^{\bra{\tau}}$. 
The groups $\opn{Aut}(\vec{\Delta})$ and
$\opn{Aut}(\vec{\mbb{Z}} \vec{\Delta})^{\bra{\tau}}$
are described in Table \tup{\ref{tab1}}. 
\end{thm}

\begin{table}
\begin{tabular}{|c|c|c|c|}
\hline
Type &  
$\opn{Aut}(\vec{\Delta})$ &
$\opn{Aut}(\vec{\mbb{Z}} \vec{\Delta})^{\bra{\tau}}$ 
& Relation 
\stick \\  \hline \hline
$A_{n}$, $n$ even & $1$ &
$\bra{\tau, \sigma}$  $\cong$  $\mbb{Z}$ & 
$\tau^{n + 1} = \sigma^{-2}$ 
\stick \\ \hline
$A_{n}$, $n$ odd & $1$ &
$\bra{\tau, \sigma} \cong \mbb{Z} \times (\mbb{Z}/2 \mbb{Z})$ 
& $\tau^{n + 1} = \sigma^{-2}$ 
\stick \\ \hline
$D_{4}$ & $S_{3}$ &
$\opn{Aut}(\vec{\Delta}) \times \bra{\tau}  
\cong S_{3} \times \mbb{Z}$ & 
$\tau^{3} = \sigma^{-1}$ 
\stick \\ \hline
$D_{n}$, $n \geq 5$ & $S_{2}$ & 
$\opn{Aut}(\vec{\Delta}) \times \bra{\tau} 
\cong S_{2} \times \mbb{Z}$
& $\tau^{n-1} = \theta \sigma^{-1}$, $n$ odd \\
& & & $\tau^{n-1} = \sigma^{-1}$, $n$ even
\stick \\ \hline
$E_{6}$ & $S_{2}$ &
$\opn{Aut}(\vec{\Delta}) \times \bra{\tau}
\cong S_{2} \times \mbb{Z}$ & 
$\tau^{6} = \theta \sigma^{-1}$ 
\stick \\ \hline
$E_{7}$ & $1$ & 
$\bra{\tau} \cong {\mbb Z}$ 
& $\tau^{9} = \sigma^{-1}$ 
\stick \\ \hline
$E_{8}$ & $1$ & 
$\bra{\tau} \cong \mbb{Z}$ 
& $\tau^{15} = \sigma^{-1}$
\stick \\ \hline
\end{tabular}

\medskip
\caption{The group 
$\opn{Aut}(\vec{\mbb{Z}} \vec{\Delta})^{\bra{\tau}}$
for a Dynkin quiver. The orientation of $\vec{\Delta}$ is shown
in Figure \ref{fig1}. In types $D_{n}$ and $E_{6}$, $\theta$ 
is the element of order $2$ in $\opn{Aut}(\vec{\Delta})$.}
\label{tab1}
\end{table}

\begin{figure}
\choosegraphics{
\[ \UseTips
\begin{xy} 
            (0,0)*+@{*}="01"*+!U{\scrp{(0,1)}},"01"
\ar@{->}    "01";"01"+(10,10)*+@{*}="02"*++!R{\scrp{(0,2)}},"02" 
\ar@{->}    "02";"02"+(10,10)*+@{*}="03"*+!D{\scrp{(0,3)}},"03" 
\ar@{->}    "02";"01"+(20,0)*+@{*}="11"*+!U{\scrp{(1,1)}},"11" 
\ar@{->}    "03";"02"+(20,0)*+@{*}="12"*++!R{\scrp{(1,2)}},"12" 
\ar@{->}    "11";"12" 
\ar@{->}    "12";"12"+(10,10)*+@{*}="13"*+!D{},"13" 
\ar@{->}    "12";"11"+(20,0)*+@{*}="21"*+!U{\scrp{(2,1)}},"21" 
\ar@{->}    "13";"12"+(20,0)*+@{*}="22"*++!R{},"22" 
\ar@{->}    "21";"22" 
\ar@{->}    "22";"22"+(10,10)*+@{*}="23"*+!D{},"23" 
\ar@{->}    "22";"21"+(20,0)*+@{*}="31"*++!L{\cdots},"31" 
\ar@{->}    "23";"22"+(20,0)*+@{*}="32"*++!L{\cdots},"32" 
\ar@{->}    "31";"32"
\ar@{->}    "32";"32"+(10,10)*+@{*}="33"*++!L{\cdots},"33" 
\ar@{->}    "02"-(20,0)*+@{*}="-12"*++!R{\cdots},"-12";"01"
\ar@{->}    "03"-(20,0)*+@{*}="-13"*++!R{\cdots},"-13";"02"
\ar@{->}    "01"-(20,0)*+@{*}="-11"*++!R{\cdots},"-11";"-12"
\ar@{->}    "-12";"-13"
\end{xy} \]
}{\includegraphics[clip]{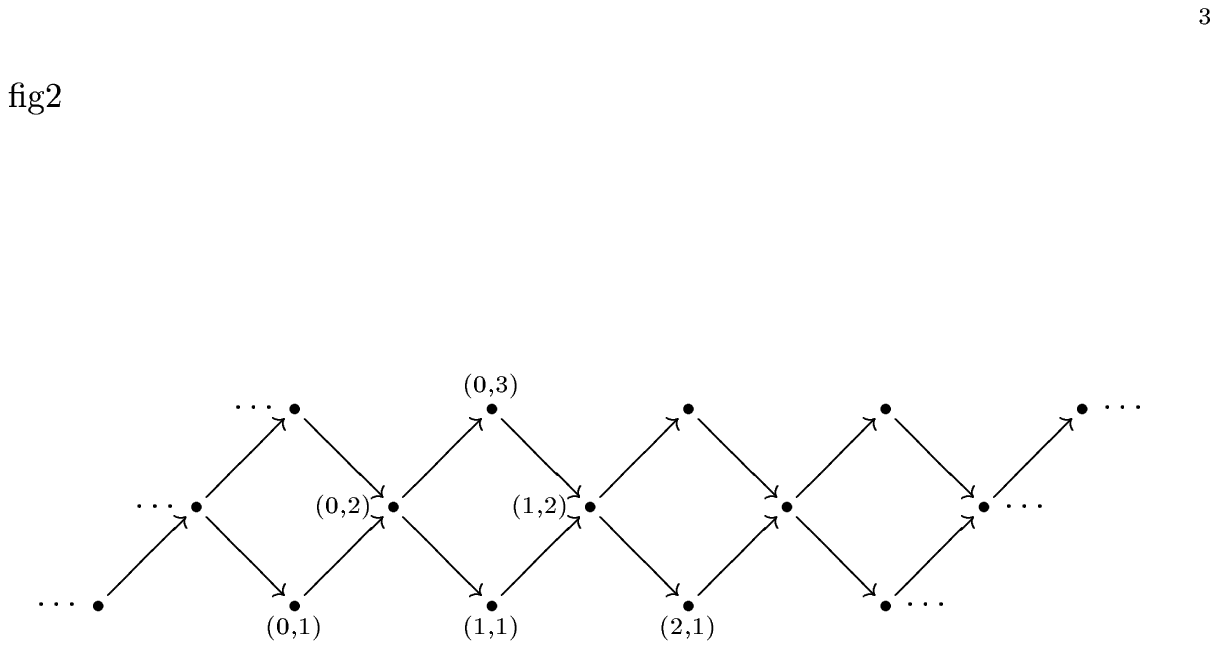}}
\caption{The quiver $\vec{\mbb{Z}} \vec{\Delta}$ for $\vec{\Delta}$
of type $A_{3}$. The vertices in 
$\cat{mod} A$ are labeled.}
\label{fig2}
\end{figure}

\begin{figure}
\choosegraphics{
\[ \UseTips
\begin{xy} 
            (0,0)*+@{*}="01"*+++!U{\scrp{(0,1)}}
            *++!R{\cdots},"01"
\ar@{->}    "01";"01"+(10,10)*+@{*}="02"*+!D{\scrp{(0,2)}},"02" 
\ar@{->}    "01";"01"+(10,0)*+@{*}="03"*+!D{\scrp{(0,3)}},"03" 
\ar@{->}    "01";"01"+(10,-10)*+@{*}="04"*+!U{\scrp{(0,4)}},"04" 
\ar@{}      (0,0)+(20,0)*+@{*}="11"*+++!U{\scrp{(1,1)}},"11"
\ar@{->}    "11";"11"+(10,10)*+@{*}="12"*+!D{\scrp{(1,2)}},"12" 
\ar@{->}    "11";"11"+(10,0)*+@{*}="13"*+!D{\scrp{(1,3)}},"13" 
\ar@{->}    "11";"11"+(10,-10)*+@{*}="14"*+!U{\scrp{(1,4)}},"14" 
\ar@{->}    "02";"11" 
\ar@{->}    "03";"11" 
\ar@{->}    "04";"11" 
\ar@{}      "11"+(20,0)*+@{*}="21"*+++!U{\scrp{(2,1)}},"21"
\ar@{->}    "21";"21"+(10,10)*+@{*}="22"*+!D{\scrp{(2,2)}},"22" 
\ar@{->}    "21";"21"+(10,0)*+@{*}="23"*+!D{\scrp{(2,3)}},"23" 
\ar@{->}    "21";"21"+(10,-10)*+@{*}="24"*+!U{\scrp{(2,4)}},"24" 
\ar@{->}    "12";"21" 
\ar@{->}    "13";"21" 
\ar@{->}    "14";"21"
\ar@{}      "21"+(20,0)*+@{*}="31"
\ar@{->}    "31";"31"+(10,10)*+@{*}="32"*++!L{\cdots},"32"
\ar@{->}    "31";"31"+(10,0)*+@{*}="33"*++!L{\cdots},"33"
\ar@{->}    "31";"31"+(10,-10)*+@{*}="34"*++!L{\cdots},"34"
\ar@{->}    "22";"31" 
\ar@{->}    "23";"31" 
\ar@{->}    "24";"31"
\end{xy} \]
}{ \includegraphics[clip]{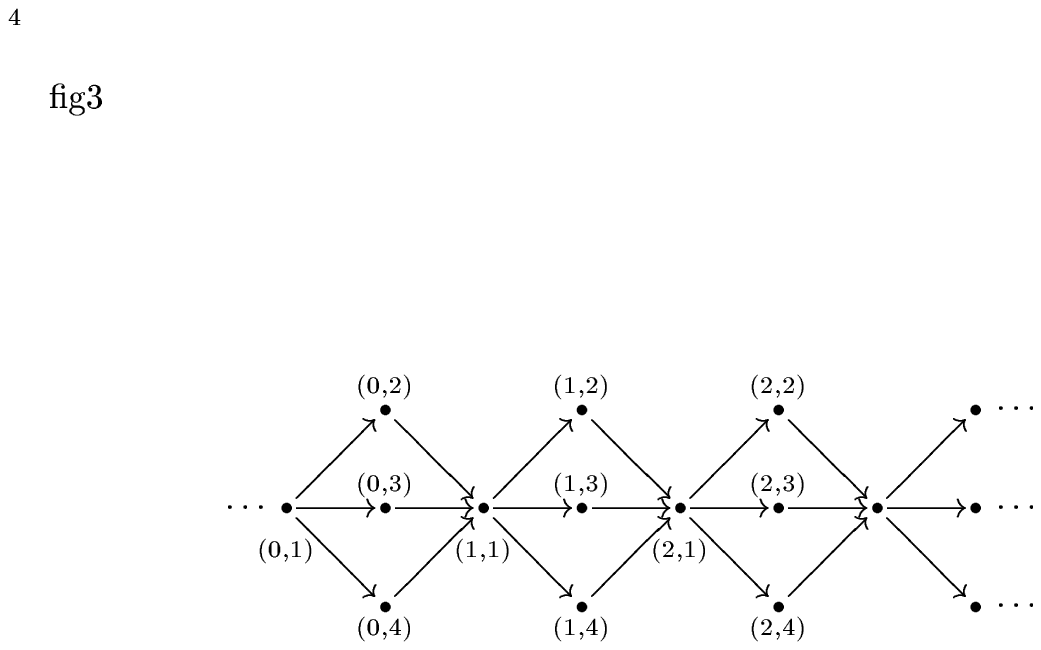}}
\caption{
The quiver $\vec{\mbb{Z}} \vec{\Delta}$ for $\vec{\Delta}$
of type $D_{4}$. The vertices in $\cat{mod} A$ are labeled.}
\label{fig3}
\end{figure}

\begin{proof}
The isomorphisms are by Theorem \ref{thm0.2} and Proposition
\ref{prop1.1}. The data in the third column of Table \ref{tab1} 
was calculated in \cite{Rn} Section 4, except for the shift $\sigma$ 
which did not appear in that paper. So we have to do a few 
calculations involving $\sigma$. Below are the
calculations for types $A_{n}$ and $D_{4}$; the rest are similar 
and are left to the reader as an exercise.

\medskip \noindent 
Type $A_{n}$: Choose the orientation in Figure \ref{fig1}. 
The quiver $\vec{\mbb{Z}} \vec{\Delta}$ looks like Figure \ref{fig2}. 
Therefore 
$\opn{Aut}(\vec{\mbb{Z}} \vec{\Delta})^{\bra{\tau}} 
= \bra{\tau, \eta}$
where $\eta(0, 1) = (0, n)$ and
$\eta(0, n) = (n - 1, 1)$. 

Now by \cite{Ha} Section I.5.5 and \cite{ARS} Sections VII.1 and
VIII.5, the quiver 
$\vec{\Gamma}(\cat{mod} A) \subset \vec{\mbb{Z}} \vec{\Delta}$
is the full subquiver on the vertices in the triangle
$\{ (m, i) \mid m \geq 0,\ m + i \leq n \}$. 
The projective vertices are
$(0, i)$ and the injective vertices are 
$(n - i, i)$, where $i \in \{ 1, \ldots, n \}$.
We see that $\sigma(0, i) = (i, n + 1 - i)$, and the quiver 
$\vec{\Gamma} \bigl( (\cat{mod} A)[1] \bigr) = 
\sigma \bigl( \vec{\Gamma}(\cat{mod}A) \bigr)$
is the full subquiver on the vertices in the triangle
$\{ (m, i) \mid m \leq n,\ m + i \geq n + 1 \}$.
Hence $\eta = \tau \sigma$ and 
$\opn{Aut}(\vec{\mbb{Z}} \vec{\Delta})^{\bra{\tau}} 
= \bra{\tau, \sigma}$.
The relation $\tau^{-(n + 1)} = \sigma^{2}$ is easily verified.

\medskip \noindent 
Type $D_{4}$: The quiver $\vec{\mbb{Z}} \vec{\Delta}$ 
is in Figure \ref{fig3}, and 
$\vec{\Gamma}(\cat{mod}A) \subset \vec{\mbb{Z}} \vec{\Delta}$ 
is a full subquiver. From the shape of $\vec{\Delta}$ we know that 
$\cat{mod} A$ should have $4$ indecomposable projective modules, 
$3$ having length $2$ and one of them simple. From the shape of the 
opposite quiver $\vec{\Delta}^{\circ}$ we also know that 
$\cat{mod} A$ should have $4$ indecomposable injective modules, $3$ 
of them simple and one of length $4$. 
Counting dimensions using Auslander-Reiten sequences we conclude 
that  $\vec{\Gamma}(\cat{mod}A)$ is the full subquiver on the 
vertices
$\{ 0, 1, 2 \} \times \vec{\Delta}_{0}$. The projective vertices 
are $\{ (0, 1), (0, i) \}$,
the injective vertices are $\{ (2, 1), (2, i) \}$, and
the simple vertices are $\{ (0, 1), (2, i) \}$,
where $i \in \{ 2, 3, 4 \}$.

For $i \in \{ 1, 2, 3, 4 \}$ let $P_{i}$, $S_{i}$ and $I_{i}$,  be 
the projective, simple and injective modules respectively, indexed 
such that $P_{i} \surj S_{i} \inj I_{i}$, and with
$P_{i} = M_{(0, i)}$. So $P_{1} = S_{1}$ and 
$I_{i} = S_{i}$ for $i  \in \{ 2, 3, 4 \}$. By the symmetry of the 
quiver it follows that there is a nonzero morphism
$(0, i) \to (2, i)$ in $k \bra{\vec{\mbb{Z}} \vec{\Delta}}$
for $i  \in \{ 2, 3, 4 \}$, and hence 
$M_{(2, i)} \cong S_{i}$

The rule for connecting $\vec{\Gamma}(\cat{mod}A)$ with 
$\vec{\Gamma}(\cat{mod}A[1])$ (see \cite{Ha} Section I.5.5)
implies that 
$M_{(3, 1)} \cong M_{(0, 1)}[1] = P_{1}[1]$. Therefore 
$M_{(3, i)} \cong P_{i'}[1]$
for $i, i' \in \{ 2, 3, 4 \}$. Now for each such $i$ there is 
an Auslander-Reiten triangle 
$M_{(2, i)} \to M_{(3, 1)} \to M_{(3, i)} \to M_{(2, i)}[1]$.
When this triangle is turned it gives an exact sequence
$0 \to P_{1} \to P_{i'} \to S_{i} \to 0$, and hence $i' = i$.
The conclusion is that $\sigma(m, i) = (m + 3, i)$ for all
$(m, i) \in (\vec{\mbb{Z}} \vec{\Delta})_{0}$, so 
$\sigma = \tau^{-3}$. 
\end{proof}

\begin{figure} 
\choosegraphics{
\[ \UseTips
\begin{array}{l}
\tilde{D}_{4} \quad
\begin{xy} 
            (0,0)*+@{*}="1"*+!R{\scrp{1}},"1"
\ar@{->}    (10,6)*+@{*}="2"*+!L{\scrp{2}},"2" 
\ar@{->}    (10,2)*+@{*}="3"*+!L{\scrp{3}},"3"
\ar@{->}    (10,-2)*+@{*}="4"*+!L{\scrp{4}},"4" 
\ar@{->}    (10,-6)*+@{*}="5"*+!L{\scrp{5}},"5" 
\end{xy} \hspace{15mm}
\begin{xy} 
            (0,0)*+@{*}="1"*+!U{\scrp{1}},"1"
\ar@{->}    (8,0)*+@{*}="2"*+!U{\scrp{2}},"2" 
\ar@{->}    "2";(16,0)*+@{*}="3"*+!U{\scrp{3}},"3" 
\ar@{->}    "3";(24,0)*+@{*}="4"*+!U{\scrp{4}},"4" 
\ar@{->}    "4";(32,0)*+@{*}="5"*+!U{\scrp{5}},"5" 
\ar@{->}    "5";(40,0)*+@{*}="6"*+!U{\scrp{6}},"6" 
\ar@{->}    "6";(48,0)*+@{*}="7"*+!U{\scrp{7}},"7" 
\ar@{->}    "7";(56,0)*+@{*}="8"*+!U{\scrp{8}},"8" 
\ar@{->}    "3";(16,8)*+@{*}="9"*+!R{\scrp{9}},"9" 
\ar@{}      (8,6)*{\tilde{E}_{8}}
\end{xy} \\
\begin{xy} 
            (0,0)*+@{*}="3"*+!R{\scrp{3}},"3"
\ar@{->}    (0,8)*+@{*}="1"*+!R{\scrp{1}},"1" 
\ar@{->}    (0,-8)*+@{*}="2"*+!R{\scrp{2}},"2"
\ar@{->}    (8,0)*+@{*}="4"*+!U{\scrp{4}};"3" 
\ar@{->}    "4";(16,0)*+@{}="5","5" 
\ar@{}      "5";(24,0)*+@{}="6","6" |*{\cdots}
\ar@{->}    (32,0)*+@{*}="2m-2";"6"
\ar@{->}    "2m-2";(40,0)*+@{*}="2m-1"*+!L{\scrp{2m-1}},"2m-1"
\ar@{->}    "2m-1";(40,8)*+@{*}="2m"*+!L{\scrp{2m}},"2m"
\ar@{->}    "2m-1";(40,-8)*+@{*}="2m+1"*+!L{\scrp{2m+1}},"2m+1"
\ar@{}      (10,6)*{\tilde{D}_{2m}}
\end{xy} \hspace{15mm} 
\begin{xy} 
            (0,0)*+@{*}="1"*+!U{\scrp{1}},"1"
\ar@{->}    (8,0)*+@{*}="2"*+!U{\scrp{2}};"1"
\ar@{->}    (16,0)*+@{*}="3"*+!U{\scrp{3}};"2" 
\ar@{->}    (24,0)*+@{*}="4"*+!U{\scrp{4}};"3" 
\ar@{->}    "4";(32,0)*+@{*}="5"*+!U{\scrp{5}},"5" 
\ar@{->}    "5";(40,0)*+@{*}="6"*+!U{\scrp{6}},"6" 
\ar@{->}    "6";(48,0)*+@{*}="7"*+!U{\scrp{7}},"7" 
\ar@{->}    "4";(24,8)*+@{*}="8"*+!R{\scrp{8}},"8" 
\ar@{}      (8,6)*{\tilde{E}_{7}}
\end{xy}  \\
\begin{xy} 
            (0,0)*+@{*}="3"*+!R{\scrp{3}},"3"
\ar@{->}    (0,8)*+@{*}="1"*+!R{\scrp{1}},"1" 
\ar@{->}    (0,-8)*+@{*}="2"*+!R{\scrp{2}},"2"
\ar@{->}    (8,0)*+@{*}="4"*+!U{\scrp{4}};"3" 
\ar@{->}    "4";(16,0)*+@{}="5","5" 
\ar@{}      "5";(24,0)*+@{}="6","6" |*{\cdots}
\ar@{->}    "6";(32,0)*+@{*}="2m-1","2m-1"
\ar@{->}    (40,0)*+@{*}="2m"*+!L{\scrp{2m}},"2m";"2m-1"
\ar@{->}    "2m";(40,8)*+@{*}="2m+1"*+!L{\scrp{2m+1}},"2m+1"
\ar@{->}    "2m";(40,-8)*+@{*}="2m+2"*+!L{\scrp{2m+2}},"2m+2"
\ar@{}      (10,6)*{\tilde{D}_{2m+1}}
\end{xy} \hspace{15mm}
\begin{xy} 
            (0,0)*+@{*}="1"*+!U{\scrp{1}},"1"
\ar@{->}    (8,0)*+@{*}="2"*+!U{\scrp{2}};"1"
\ar@{->}    (16,0)*+@{*}="3"*+!U{\scrp{3}};"2" 
\ar@{->}    "3";(24,0)*+@{*}="4"*+!U{\scrp{4}},"4" 
\ar@{->}    "4";(32,0)*+@{*}="5"*+!U{\scrp{5}},"5" 
\ar@{->}    "3";(16,8)*+@{*}="6"*+!R{\scrp{6}},"6" 
\ar@{->}    "6";(16,16)*+@{*}="7"*+!R{\scrp{7}},"7" 
\ar@{}      (8,6)*{\tilde{E}_{6}}
\end{xy} 
\end{array} \]
}{\includegraphics[clip]{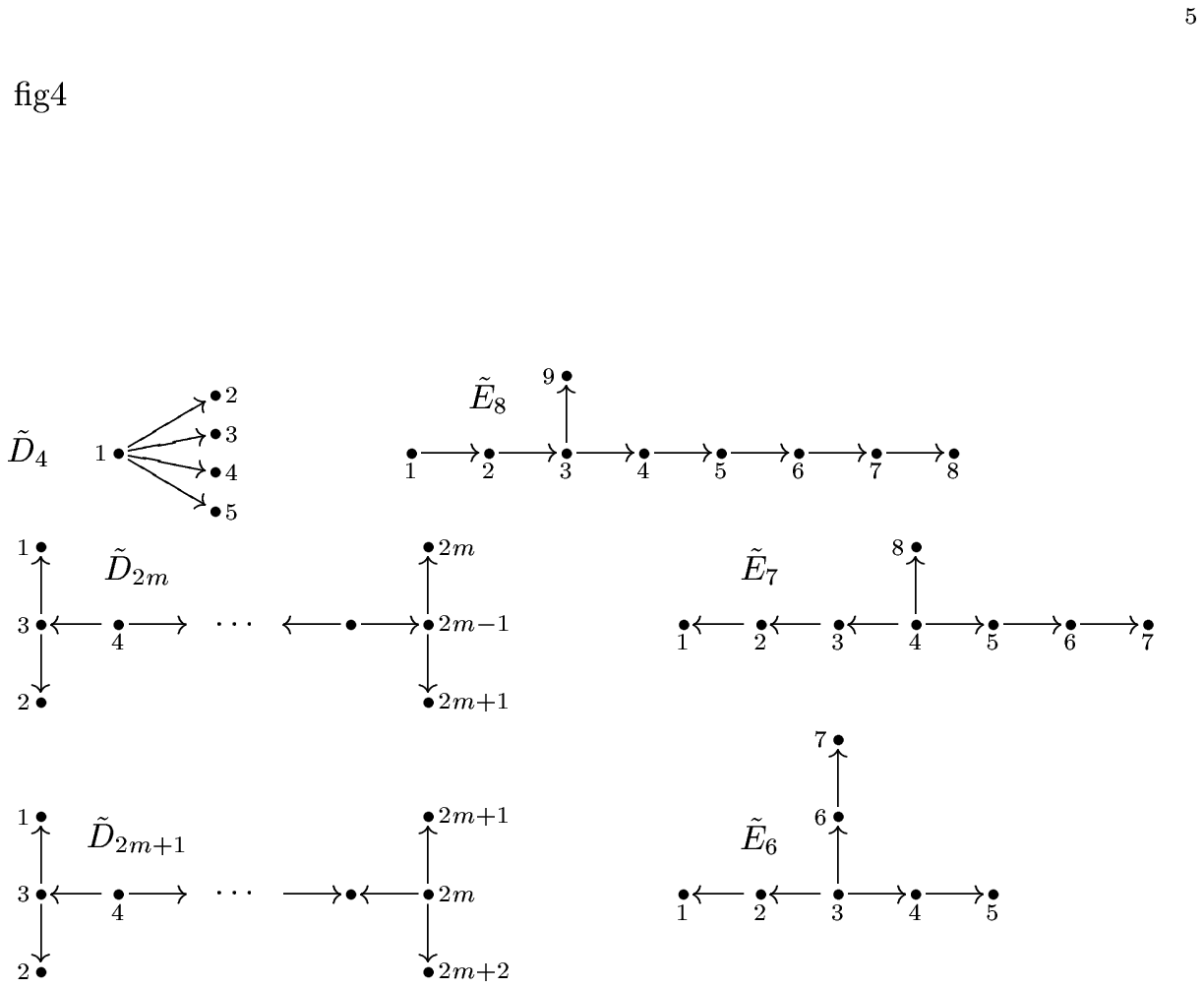}}
\caption{Orientations for the affine tree graphs} \label{fig4}
\end{figure}

\begin{thm} \label{thm4.2}
Let $\vec{\Delta}$ be a quiver of type $\tilde{D}_{n}$, 
$\tilde{E}_{6}$, $\tilde{E}_{7}$ or $\tilde{E}_{8}$, with the 
orientation shown in Figure \tup{\ref{fig4}}. Then 
$\opn{Pic}_{k}(A) \cong \opn{Aut}(\vec{\Delta})$
and
\[ \opn{DPic}_{k}(A) \cong \mbb{Z} \times
\opn{Aut}(\vec{\mbb{Z}} \vec{\Delta})^{\bra{\tau}} . \]
The structure of the group
$\opn{Aut}(\vec{\mbb{Z}} \vec{\Delta})^{\bra{\tau}}$ 
is given in Table \tup{\ref{tab2}}.
\end{thm}

\begin{proof}
The isomorphisms follow from Theorem \ref{thm0.2} and Proposition 
\ref{prop1.1}. The structure of 
$\opn{Aut}(\vec{\mbb{Z}} \vec{\Delta})^{\bra{\tau}}$
is quite easy to check in all cases. In type $\tilde{D}_{n}$, 
$n \geq 5$ odd, the automorphism 
$\eta \in \opn{Aut}(\vec{\mbb{Z}} \vec{\Delta})^{\bra{\tau}}$ 
is 
\[
\eta(i,j) = 
\begin{cases}
(i, n+2 - j) & \text{if } j = 2, n \\
(i - \frac{1-(-1)^j}{2}, n+2 - j) & \text{otherwise} .
\end{cases}
\]
\end{proof}

\begin{table}
\begin{tabular}{|c|c|c|c|}
\hline
Type &  $\opn{Aut}(\vec{\Delta})$ &
$\opn{Aut}(\vec{\mbb{Z}} \vec{\Delta})^{\bra{\tau}}$ &
Relations \stick \\ \hline \hline
$\tilde{D}_{4}$ & $S_{4}$ & 
$\opn{Aut}(\vec{\Delta}) \times \bra{\tau} \cong
S_{4} \times \mbb{Z}$ &
\stick \\ \hline
$\tilde{D}_{n}$, $n \geq 5$ even &
$S_{2} \ltimes S_{2}^{2}$ &
$\opn{Aut}(\vec{\Delta}) \times \bra{\tau} \cong 
(S_{2} \ltimes S_{2}^{2}) \times \mbb{Z}$ &
\stick \\ \hline
$\tilde{D}_{n}$, $n \geq 5$ odd &
$S_{2}^{2}$ &
$\opn{Aut}(\vec{\Delta}) \times \bra{\eta} \cong 
S_{2}^{2} \times \mbb{Z}$ & $\eta^{2} = \tau$
\stick \\ \hline
$\tilde{E}_{6}$ & $S_{3}$ &
$\opn{Aut}(\vec{\Delta}) \times \bra{\tau} \cong 
S_{3} \times \mbb{Z}$ &
\stick \\ \hline
$\tilde{E}_{7}$ & $S_{2}$ &
$\opn{Aut}(\vec{\Delta}) \times \bra{\tau} \cong 
S_{2} \times \mbb{Z}$ &
\stick \\ \hline
$\tilde{E}_{8}$ & $1$ & 
$\bra{\tau} \cong \mbb{Z}$ &
\stick \\ \hline
\end{tabular}

\medskip
\caption{The groups
$\opn{Aut}(\vec{\mbb{Z}} \vec{\Delta})^{\bra{\tau}}$ 
for the affine tree quivers shown in Figure \tup{\ref{fig4}}.}
\label{tab2}
\end{table}

\begin{figure} 
\choosegraphics{
\[ \UseTips
\vec{\Omega}_{n} \quad 
\begin{xy} 
    (0,0)*+@{*}="1"*+!CR{\scrp{1}},"1"
\ar@/^/@<2.5ex>@{->}|*+{\scrp{\alpha_{1}}}     
    (20,0)*+@{*}="2"*+!CL{\scrp{2}},"2"
\ar@/^/@<1ex>@{->}|*+{\scrp{\alpha_{2}}}     "2"
\ar@{}|{\vdots}     "2"
\ar@/_/@<-2ex>@{->}|*+{\scrp{\alpha_{n}}}    "2"
\end{xy} 
\hspace{25mm} 
\vec{T}_{p, q} \quad
\begin{xy} 
    (0,10)*+@{*}="1"*+!DC{\scrp{1}},"1"
\ar@{->}^{\alpha_{1}}  (10,10)*+@{*}="2"*+!DC{\scrp{2}},"2"
\ar@{->}^{\alpha_{2}}  "2";(20,10)="3"
\ar@{}|*{\cdots}       "3";(30,10)="4"
\ar@{->}^{\alpha_{p - 1}}  "4";(40,10)*+@{*}="p"*+!DC{\scrp{p}}
\ar@{->}^{\alpha_{p}}  
    "p";(40,0)*+@{*}="p+1"*+!UC{\scrp{p+1}},"p+1"
\ar@{->}_{\beta_{q}}  
    "1";(0,0)*+@{*}="p+q"*+!UC{\scrp{p+q}},"p+q"
\ar@{->}^{\beta_{q-1}}  
    "p+q";(10,0)*+@{*}="p+q-1"*+!UC{\scrp{p+q-1}},"p+q-1"
\ar@{->}^{\beta_{q-2}}  "p+q-1";(20,0)="p+q-2"
\ar@{}|*{\cdots}       "p+q-2";(30,0)="p+2"
\ar@{->}^{\beta_{1}}  "p+2";"p+1"
\end{xy} \]
}{\includegraphics[clip]{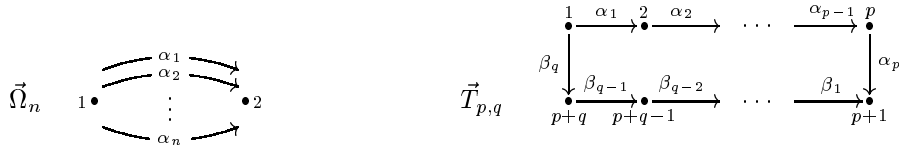}}
\caption{The quivers $\vec{\Omega}_{n}$ and $\vec{T}_{p, q}$.} 
\label{fig5}
\end{figure}

\begin{thm} \label{thm4.3}
For any $n \geq 2$ let $\vec{\Omega}_{n}$ be the quiver shown in 
Figure \tup{\ref{fig5}}, and let $A := k \vec{\Omega}_{n}$ be the 
path algebra. Then 
$\opn{Pic}_{k}(A) \cong \opn{PGL}_{n}(k)$
and
\[ \opn{DPic}_{k}(A) \cong \mbb{Z} \times \bigl( \mbb{Z} \ltimes
\opn{PGL}_{n}(k) \bigr) . \]
In the semidirect product the action of a generator 
$\rho \in \mbb{Z}$ on a matrix $F \in \opn{PGL}_{n}(k)$ is 
$\rho F \rho^{-1} = (F^{-1})^{\mrm{t}}$. 
\end{thm}

\begin{proof}
As in the proof of Lemma \ref{lem1.1} and Proposition \ref{prop1.1},
the group of auto-equivalences of the path category is 
$\opn{Aut}_{k}(k \bra{\vec{\Omega}_{n}}) = 
\opn{Aut}_{k}^{0}(k \bra{\vec{\Omega}_{n}}) \cong 
\opn{GL}_{n}(k)$.
Hence
$\opn{Pic}_{k}(A) \cong 
\opn{Out}_{k}(k \bra{\vec{\Omega}_{n}}) \cong \opn{PGL}_{n}(k)$.

Given $F \in \opn{Aut}_{k}(k \bra{\vec{\Omega}_{n}})$
let 
$[a_{i, j}] \in \opn{GL}_{n}(k)$ be its matrix w.r.t.\ to the basis
$\{ \alpha_{i} \}$, and let 
$[b_{i, j}] := ([a_{i, j}]^{-1})^{\mrm{t}}$. 
Define an auto-equivalence 
$\tilde{F} \in 
\opn{Aut}_{k}^{0}(k \bra{\vec{\mbb{Z}} \vec{\Omega}_{n}})$
with
$\tilde{F}(m, \alpha_{i}) = \sum_{j} a_{i, j} (m, \alpha_{j})$
and
$\tilde{F}(m, \alpha^{*}_{i}) = 
\sum_{j} b_{i, j} (m, \alpha^{*}_{j})$, $m \in \mbb{Z}$.
Then $\tilde{F}$ preserves all mesh relations, and by a linear algebra 
argument we see that up to scalars at each vertex, the only 
elements of 
$\opn{Aut}^{0}_{k}(k \bra{\vec{\mbb{Z}} \vec{\Omega}_{n}})$ 
are of the form $\tilde{F}$. 

Let $\rho \in \opn{Aut}(\vec{\mbb{Z}} \vec{\Omega}_{n})$ be 
$\rho(m, 1) = (m, 2)$ and $\rho(m, 2) = (m + 1, 1)$, with the 
obvious action on arrows to make it commute with the polarization 
$\mu$. Then 
$\opn{Out}_{k}(k \bra{\vec{\mbb{Z}} \vec{\Omega}_{n}, I_{\mrm{m}}})$
is generated by $\opn{PGL}_{n}(k)$ and $\rho$, so
$\opn{Out}_{k}(k \bra{\vec{\mbb{Z}} \vec{\Omega}_{n}, I_{\mrm{m}}}) 
\cong \mbb{Z} \ltimes \opn{PGL}_{n}(k)$.
The formula for $\tilde{F}$ above shows that 
$\rho F \rho^{-1} = (F^{-1})^{\mrm{t}}$ for 
$F \in \opn{PGL}_{n}(k)$.

Finally use Theorem \ref{thm3.1}.
\end{proof}

\begin{rem}
By \cite{Be} and \cite{BO} we see that for $n = 2$ in the theorem 
above, 
$\opn{DPic}_{k}(A) \cong \mbb{Z} \times \mbb{Z} \times
\opn{PGL}_{2}(k)$. 
The apparent discrepancy is explained by the fact that 
$\mbb{Z} \ltimes \opn{PGL}_{2}(k) \cong
\mbb{Z} \times \opn{PGL}_{2}(k)$
via
$(m, F) \mapsto (m, H^{m} F)$,
where
$H = \left[ \begin{smallmatrix}
0 & -1 \\ 1 & 0 
\end{smallmatrix} \right]$. 
\end{rem}

For integers $p \geq q \geq 1$ let $\vec{T}_{p, q}$ be the quiver 
shown in Figure \ref{fig5}. 
Let $\vec{\Delta}$ be a quiver with underlying graph $\tilde{A}_{n}$.
Then $\vec{\Delta}$ can be brought to one 
of the quivers $\vec{T}_{p, q}$, $p + q = n+1$, by a sequence of 
admissible reflections $s_{x}^{-}$ at source vertices (see Section
6). Therefore 
\[ \opn{DPic}_{k}(k \vec{\Delta}) \cong 
\opn{DPic}_{k}(k \vec{T}_{p, q}) . \]

\begin{thm} \label{thm4.4}
Let $A$ be the path algebra $k \vec{T}_{p, q}$. 
\begin{enumerate}
\item If $p = q = 1$ then 
$\opn{Pic}_{k}(A) \cong \opn{PGL}_{2}(k)$
and
$\opn{DPic}_{k}(A) \cong \mbb{Z} \times 
\bigl( \mbb{Z} \ltimes \opn{PGL}_{2}(k) \bigr)$.
\item If $p > q = 1$ then 
$\opn{Pic}_{k}(A) \cong
\left[ \begin{smallmatrix}
k^{\times} & k \\ 0 & 1
\end{smallmatrix} \right]$
and
$\opn{DPic}_{k}(A) \cong \mbb{Z} \times 
\bigl( \mbb{Z} \ltimes 
\left[ \begin{smallmatrix} k^{\times} & k \\ 0 & 1 
\end{smallmatrix} \right] \bigr)$. 
\item If $p = q > 1$ then 
$\opn{Pic}_{k}(A) \cong S_{2} \ltimes k^{\times}$
and
$\opn{DPic}_{k}(A) \cong \mbb{Z}^{2} \times 
\bigl( S_{2} \ltimes k^{\times} \bigr)$.
\item If $p > q > 1$ then 
$\opn{Pic}_{k}(A) \cong k^{\times}$
and
$\opn{DPic}_{k}(A) \cong \mbb{Z}^{2} \times 
k^{\times}$.
\end{enumerate}
\end{thm}

\begin{proof}
1. This is because 
$\vec{T}_{1, 1} = \vec{\Omega}_{2}$. 

\medskip \noindent 2.
Here the group of auto-equivalences of $k \bra{\vec{T}_{p, q}}$
is, in the notation of the proof of Proposition \ref{prop1.1},
$\opn{Aut}_{k}(k \bra{\vec{T}_{p, q}}) \cong 
(k^{\times})^{p + 1} \times k$,
and the group of isomorphisms is $(k^{\times})^{p}$. Therefore
$\opn{Out}_{k}(k \bra{\vec{T}_{p, q}})$
is isomorphic to $k^{\times} \times k$ as varieties, and as matrix 
group
$\opn{Out}_{k}(k \bra{\vec{T}_{p, q}})
\cong \left[ \begin{smallmatrix} k^{\times} & k \\ 0 & 1 
\end{smallmatrix} \right]$.
The auto-equivalence associated to 
$\left[ \begin{smallmatrix} a & b \\ 0 & 1 
\end{smallmatrix} \right] \in 
\left[ \begin{smallmatrix} k^{\times} & k \\ 0 & 1 
\end{smallmatrix} \right]$
is $\alpha_{i} \mapsto \alpha_{i}$ and
$\beta_{1} \mapsto a \beta_{1} + b \alpha_{p} \cdots \alpha_{1}$.

The quiver $\vec{\mbb{Z}} \vec{T}_{p, q}$ has no multiple arrows. Let 
$\rho$ be the symmetry $\rho(m, i) = (m, i - 1)$ for $i \geq 2$, 
and $\rho(m, 1) = (m - 1, p)$. Then $\rho$ generates 
$\opn{Aut}(\vec{\mbb{Z}} \vec{T}_{p, q})^{\bra{\tau}}$, and we can use
Theorem \ref{thm0.2}. The action of $\rho$ on 
$\opn{Out}_{k}(k \bra{\vec{T}_{p, q}})$ is
$\rho 
\left[ \begin{smallmatrix} a & b \\ 0 & 1 \end{smallmatrix} \right]
\rho^{-1} = 
\left[ \begin{smallmatrix} a & -b \\ 0 & 1 \end{smallmatrix} 
\right]$.

\medskip \noindent 3.
Here  
$\opn{Aut}^{0}_{k}(k \bra{\vec{T}_{p, q}}) \cong 
(k^{\times})^{2p}$,
and the subgroup of isomorphisms is $(k^{\times})^{2p - 1}$. 
The symmetry $\theta \in \opn{Aut}(\vec{T}_{p, q})$ of order $2$ 
acts on $k^{\times}$ by 
$\theta a \theta^{-1} = a^{-1}$. 

Let $\rho$ be the symmetry  
$\rho(m, 1) = (m - 1, p + q)$,
$\rho(m, i) = (m, i - 1)$ if $2 \leq i \leq p + 1$,
and
$\rho(m, i) = (m - 1, i - 1)$ if $p + 2 \leq i \leq p + q$.
Then $\rho$ and $\theta$ commute, and they generate 
$\opn{Aut}(\vec{\mbb{Z}} \vec{T}_{p, q})^{\bra{\tau}}$. The
action of $\rho$ on $\opn{Aut}_{k}(k \bra{\vec{T}_{p, q}})$ 
is trivial. 
 
\medskip \noindent 4.
Similar to case 3. 
\end{proof}

\section{The Reflection Groupoid of a Graph}

In this section we interpret the reflection functors of \cite{BGP} 
and the tilting modules of \cite{APR} in the setup of derived 
categories. 

Let $\Delta$ be a tree with $n$ vertices. Denote by 
$\opn{Or}(\Delta)$ the set of orientations of the edge set 
$\Delta_{1}$. For $\omega \in \opn{Or}(\Delta)$ let 
$\vec{\Delta}_{\omega}$ be the resulting quiver, and let 
$A_{\omega}$ be the path algebra $k \vec{\Delta}_{\omega}$. 

Given two orientations $\omega,\omega'$ let
\[ \opn{DPic}_{k}(\omega,\omega') := 
\frac{ \{ \text{two-sided tilting complexes } 
T \in \cat{D}^{\mrm{b}}(\cat{Mod}(A_{\omega'} \otimes_{k}
A_{\omega}^{\circ})) \} }
{\text{isomorphism}} . \]
The {\em derived Picard groupoid} of $\Delta$ is the groupoid 
$\opn{DPic}_{k}(\Delta)$ with object set $\opn{Or}(\Delta)$ and 
morphism sets $\opn{DPic}_{k}(\omega,\omega')$. Thus when 
$\omega = \omega'$ we recover the derived Picard group 
$\opn{DPic}_{k}(A_{\omega})$.

For an orientation $\omega$ and a vertex $x$ let
$P_{x, \omega} \in \cat{mod} A_{\omega}$
be the corresponding indecomposable projective module. 
Denote by $\tau_{\omega}$ the translation functor of 
$\cat{D}^{\mrm{b}}(\cat{mod} A_{\omega})$, 
i.e.\ the functor
$\tau_{\omega} = A_{\omega}^{*}[-1] \otimes^{\mrm{L}}_{A_{\omega}}
-$.

Suppose $x \in (\vec{\Delta}_{\omega})_{0}$ is a source. Define 
$s^{-}_{x} \omega$ to be the orientation obtained from $\omega$ by 
reversing the arrows starting at $x$. Let
\[ T_{x, \omega} := \tau_{\omega}^{-1} P_{x, \omega} \oplus
\bigl( \bigoplus_{y \neq x} P_{y, \omega} \bigr) 
\in \cat{mod} A_{\omega} . \]
According to \cite{APR} Section 3,
$T_{x, \omega}$ is a tilting module, with 
$\opn{End}_{A_{\omega}}(T_{x, \omega})^{\circ} \cong 
A_{s^{-}_{x} \omega}$. It is called an {\em APR tilting module}. 
One has isomorphisms in $\cat{mod} A_{s^{-}_{x} \omega}$:
\begin{equation} \label{eqn5.1}
\begin{aligned}
\opn{Hom}_{A_{\omega}}(T_{x, \omega}, P_{y, \omega}) & \cong
P_{y, s^{-}_{x} \omega} \quad \text{if } y \neq x , \\
\opn{Hom}_{A_{\omega}}(T_{x, \omega}, 
\tau_{\omega}^{-1} P_{x, \omega}) & \cong
P_{x, s^{-}_{x} \omega} . 
\end{aligned}
\end{equation}
Under the anti-equivalence between $\cat{mod} A_{\omega}$ and 
the category of finite dimensional representations of the quiver
$\vec{\Delta}_{\omega}$, the reflection functor of \cite{BGP} is 
sent to
$\opn{Hom}_{A_{\omega}}(T_{x, \omega}, -) : \cat{mod} A_{\omega}
\to \cat{mod} A_{s^{-}_{x} \omega}$. 

\begin{dfn} \label{dfn5.1}
The {\em reflection groupoid} of $\Delta$ is the subgroupoid 
$\opn{Ref}(\Delta) \subset$ \linebreak
$\opn{DPic}_{k}(\Delta)$
generated by the two-sided tilting complexes 
$T_{x, \omega} \in 
\cat{D}^{\mrm{b}}(\cat{Mod} 
(A_{\omega} \otimes_{k} A_{s^{-}_{x} \omega}^{\circ}))$,
as $\omega$ runs over $\opn{Or}(\Delta)$ and $x$ runs over the 
sources in $\vec{\Delta}_{\omega}$.
\end{dfn}

Given an orientation $\omega$ the set 
$\{ [P_{x, \omega}] \}_{x \in \Delta_{0}}$
is a basis of the Grothendieck group 
$\opn{K}_{0}(A_{\omega}) = 
\opn{K}_{0}(\cat{D}^{\mrm{b}}(\cat{mod} A_{\omega}))$. 
Let $\mbb{Z}^{\Delta_{0}}$ be the free abelian group with basis 
$\{ e_{x} \}_{x \in \Delta_{0}}$. Then 
$[P_{x, \omega}] \mapsto e_{x}$ determines a canonical isomorphism
$\opn{K}_{0}(A_{\omega}) \iso \mbb{Z}^{\Delta_{0}}$. 
For a two-sided tilting complex
$T \in \opn{DPic}_{k}(\omega,\omega')$
let 
$\chi_{0}(T) : \opn{K}_{0}(A_{\omega}) \iso
\opn{K}_{0}(A_{\omega'})$
be 
$\chi_{0}(T)([M]) := [T \otimes^{\mrm{L}}_{A_{\omega}} M]$.
Using the projective bases we get a functor (when we consider a 
group as a groupoid with a single object)
\[ \chi_{0} : \opn{DPic}_{k}(\Delta) \to
\opn{Aut}_{\mbb{Z}}(\mbb{Z}^{\Delta_{0}}) \cong 
\opn{GL}_{n}(\mbb{Z}) . \]

Recall that for a vertex $x \in \Delta_{0}$ one defines the 
reflection
$s_{x} \in \opn{Aut}_{\mbb{Z}}(\mbb{Z}^{\Delta_{0}})$
by 
\[ \begin{aligned}
s_{x} e_{x} & := -e_{x} + \sum_{ \{x, y\} \in  \Delta_{1}} e_{y}
, \\
s_{x} e_{y} & := e_{y} \quad \text{if } y \neq x .  
\end{aligned} \]
The {\em Weyl group} of $\Delta$ is the subgroup  
$W(\Delta) \subset \opn{Aut}_{\mbb{Z}}(\mbb{Z}^{\Delta_{0}})$
generated by the reflections $s_{x}$. 

\begin{prop} \label{prop5.1}
Let $x$ be a source in the quiver $\vec{\Delta}_{\omega}$. Then 
\[ \chi_{0}(T_{x, \omega}) = s_{x} . \]
\end{prop}

\begin{proof}
There is an Auslander-Reiten sequence
\[ 0 \to P_{x, \omega} \to \bigoplus_{(x \to y) \in 
(\vec{\Delta}_{\omega})_{1}} P_{y, \omega} \to
\tau_{\omega}^{-1} P_{x, \omega} \to 0 \]
in $\cat{mod} A_{\omega}$. Applying the functor
$T_{x, \omega}^{\vee} \otimes^{\mrm{L}}_{A_{\omega}} - \cong
\opn{R} \opn{Hom}_{A_{\omega}}(T_{x, \omega}, -)$
to this sequence, and using formula (\ref{eqn5.1}), we get a 
triangle
\[ T_{x, \omega}^{\vee} \otimes^{\mrm{L}}_{A_{\omega}} P_{x, \omega}
\to \bigoplus_{ \{x, y\} \in \Delta_{1}} 
P_{y, s_{x}^{-} \omega}
\to P_{x, s_{x}^{-} \omega}
\to (T_{x, \omega}^{\vee} \otimes^{\mrm{L}}_{A_{\omega}}
P_{x, \omega})[1] \] 
in $\cat{D}^{\mrm{b}}(\cat{mod} A_{s_{x}^{-} \omega})$.
Hence
\[ [T_{x, \omega}^{\vee} 
\otimes^{\mrm{L}}_{A_{\omega}} P_{x, \omega}] = 
-[P_{x, s_{x}^{-} \omega}] + 
\sum_{ \{x, y\} \in \Delta_{1}} [P_{y, s_{x}^{-} \omega}]
\in \opn{K}_{0}(A_{s_{x}^{-} \omega}) . \]
On the other hand for $y \neq x$ we have
$[T_{x, \omega}^{\vee} 
\otimes^{\mrm{L}}_{A_{\omega}} P_{y, \omega}] = 
[P_{y, s_{x}^{-} \omega}]$.
This proves that
$\chi_{0}(T_{x, \omega}^{\vee}) = s_{x}$; 
but $s_{x} = s_{x}^{-1}$. 
\end{proof}

An immediate consequence is:

\begin{cor}
$\chi_{0}(\opn{Ref}(\Delta)) = W(\Delta)$.
\end{cor}

An ordering $(x_{1}, \ldots, x_{n})$ of $\Delta_{0}$ 
is called source-admissible for an orientation $\omega$ if 
$x_{i}$ is a source in the quiver
$\vec{\Delta}_{s_{x_{i - 1}}^{-} \cdots s_{x_{1}}^{-} \omega}$
for all $1 \leq i \leq n$. Any orientation has source-admissible 
orderings of the vertices. 

\begin{prop} \label{prop5.2}
Let $(x_{1}, \ldots, x_{n})$ be a source-admissible ordering of 
$\Delta_{0}$ for an orientation $\omega$. Write
$\omega_{i} := s_{x_{i}}^{-} \cdots s_{x_{1}}^{-} \omega$,
$A_i := A_{\omega_{i}}$ and 
$T_{i} := T_{x_{i}, \omega_{i - 1}}$. Then
\[ T^{\vee}_{n} \otimes^{\mrm{L}}_{A_{n - 1}} \cdots
\otimes^{\mrm{L}}_{A_{2}} T^{\vee}_{2} 
\otimes^{\mrm{L}}_{A_{1}} T^{\vee}_{1} \cong
A_{\omega}^{*}[-1] \]
in $\cat{D}^{\mrm{b}}(\cat{Mod} A_{\omega}^{\mrm{e}})$.
\end{prop}     
               
\begin{proof} 
For an orientation $\omega$ let
$\vec{\Gamma}^{\mrm{irr}}_{\omega} \subset 
\vec{\Gamma}(\cat{D}^{\mrm{b}}(\cat{mod} A_{\omega}))$
be the quiver of definition \ref{dfn2.1}. As usual 
$(\vec{\Gamma}^{\mrm{irr}}_{\omega})_{0}$ denotes the set of 
vertices of $\vec{\Gamma}^{\mrm{irr}}_{\omega}$. 
Let $G(\Delta)$ be the groupoid with object set 
$\opn{Or}(\Delta)$, and morphism sets 
$\opn{Iso}\bigl( (\vec{\Gamma}^{\mrm{irr}}_{\omega})_{0},
(\vec{\Gamma}^{\mrm{irr}}_{\omega'})_{0} \bigr)$
for $\omega, \omega' \in \opn{Or}(\Delta)$. The groupoid $G(\Delta)$ 
acts faithfully on the family of sets
$X(\Delta) := \{ 
(\vec{\Gamma}^{\mrm{irr}}_{\omega})_{0} \}_{\omega \in
\opn{Or}(\Delta)}$.
According to Theorem \ref{thm0.2} there is 
an injective map of groupoids 
$q : \opn{DPic}_{k}(\Delta) \inj G(\Delta)$. 

Let us first assume $\Delta$ is a Dynkin graph. Then there is a 
canonical isomorphism of sets
$X(\Delta) \cong \mbb{Z} \times {\Delta}_{0} \times \opn{Or}(\Delta)$.
The action of $q(\tau_{\omega})$ on $X(\Delta)$ is
$q(\tau_{\omega})(i, x, \omega) = (i - 1, x, \omega)$. By formula 
(\ref{eqn5.1}), the action of $q(T^{\vee}_{x, \omega})$ on 
$X(\Delta)$ is
$q(T^{\vee}_{x, \omega})(0, y, \omega) = 
(0, y, s_{x}^{-} \omega)$ 
if $y \neq x$, and 
$q(T^{\vee}_{x, \omega})(1, x, \omega) 
= (0, x, s_{x}^{-} \omega)$.
Since $q(\tau_{\omega})$ commutes with $q(T^{\vee}_{x, \omega})$ 
we have 
\[ q(T^{\vee}_{n} \otimes^{\mrm{L}}_{A_{n - 1}} \cdots
\otimes^{\mrm{L}}_{A_{1}} T^{\vee}_{1})(i, x, \omega) = 
(i - 1, x, \omega) = q(\tau_{\omega})(i, x, \omega)  \]
for any $x \in \Delta_{0}$ and $i \in \mbb{Z}$.

If $\Delta$ is not Dynkin then 
$X(\Delta) \cong \mbb{Z} \times \mbb{Z} \times 
{\Delta}_{0} \times \opn{Or}(\Delta)$,
$q(\tau_{\omega})(j, i, x, \omega) = (j, i - 1, x, \omega)$, 
etc., and the proof is the same after these modifications.
\end{proof} 
 
\begin{prop} \label{prop5.3}
For any orientation $\omega$,
\[ \opn{Ref}(\Delta)(\omega, \omega) = \bra{\tau_{\omega}} . \]
\end{prop} 

\begin{proof}
We will only treat the Dynkin case; the general case is proved 
similarly with modifications like in the previous proof. 

Let $T \in \opn{Ref}(\Delta)(\omega, \omega)$. From the proof 
above we see that 
$q(T)(0, x, \omega) = (i(x), x, \omega)$
for some $i(x) \in \mbb{Z}$. A quiver map
$\pi : \vec{\Delta}_{\omega} \to \vec{\mbb{Z}} \vec{\Delta}_{\omega}$
with $\pi(x) = (i(x), x)$ must have $i(x) = i$ for all $x$, 
since $\Delta$ is a tree. Therefore $q(T) = 
q(\tau_{\omega}^{-i})$.
\end{proof} 
 
\begin{rem}
The explicit calculations in Section 4 show that the shift 
$\sigma = A[1]$ is not in 
$\bra{\tau} \subset \opn{DPic}_{k}(A)$
for most algebras $A$. Thus 
$\opn{Ref}(\Delta) \subsetneqq \opn{DPic}_{k}(\Delta)$
for most graphs $\Delta$.
\end{rem}



\begin{thebibliography}{SGA6}
\bibitem[APR]{APR} M.\ Auslander, M.\ I.\ Platzeck and I.\ Reiten, 
        Coxeter functors without diagrams,
        Trans.\ AMS \textbf{250} (1979), 1-46.
\bibitem[ARS]{ARS} M.\ Auslander, I.\ Reiten and Sverre O.\ Smal{\o},
        ``Representation Theory of Artin Algebras,''
        Cambridge Studies in Advanced Math.\ \textbf{36},
        Cambridge, UK 1997 (corrected paperback edition).
\bibitem[Be]{Be} A.A. Beilinson, Coherent sheaves on $\mbf{P}^{n}$
        and problems of linear algebra, Func.\ Anal.\ Appl.\
        \textbf{12} (1978), 214-216.
\bibitem[BGP]{BGP} I.\ N.\ Bernstein, I.\ M.\ Gelfand and V.\ A.\ 
        Ponomarev, Coxeter functors and Gabriel's theorem, 
        Usp.\ Mat.\ Nauk  \textbf{28} (1973) 19-23, 
        Transl.\ Russ.\ Math.\ Surv.\ \textbf{28} (1973), 17-32.
\bibitem[BK]{BK} A.I. Bondal and M.M Kapranov, Representable 
        functors, Serre functors, and mutations, Izv.\ Akad.\ Nauk.\ 
        SSSR Ser.\ Mat.\ \textbf{53} (1989), 1183-1205; English 
        trans.\ Math.\ USSR Izv.\ \textbf{35} (1990), 519-541.  
\bibitem[BO]{BO} A.I.\ Bondal and D.O.\ Orlov,
        Reconstruction of a variety from the derived category and
        groups of autoequivalences, preprint;
        eprint: alg-geom/9712029.
\bibitem[GR]{GR} P. Gabriel and A.V. Roiter, ``Representations of 
        finite-dimensional algebras,''  
        Springer-Verlag, Berlin, 1997.
\bibitem[GS]{GS} F. Guil-Asensio and M. Saorin, 
        The automorphism group and the Picard group of a monoidal 
        algebra, Comm.\ Algebra \textbf{27} (1999), 857-887.
\bibitem[Ha]{Ha} D.\ Happel, ``Triangulated Categories in the
        Representation Theory of Finite-Dimensional Algebras,''
        London Math.\ Soc.\ Lecture Notes \textbf{119},
        University Press, Cambridge, 1987.
\bibitem[HR]{HR} D.\ Happel and C.\ M.\ Ringel, Tilted Algebras,
        Trans.\ AMS \textbf{274} (1982), 399-443.
\bibitem[Ko]{Ko} M. Kontsevich, Cours \`{a} l'\'{E}cole Normale 
        Sup\'{e}rieure, Paris, 1998.
\bibitem[KR]{KR} M. Kontsevich and A. Rosenberg, Noncommutative
        smooth spaces, preprint; eprint math.AG/9812158.
\bibitem[LM]{LM} H. Lenzing and H. Meltzer, The automorphism group of 
        the derived category for a weighted projective line,
        to appear in Comm.\ Algebra.
\bibitem[MR]{MR} J.C.\ McConnell and J.C.\ Robson, ``Noncommutative
        Noetherian Rings,'' Wiley, Chichester, 1987.
\bibitem[Or]{Or} D.O. Orlov, Equivalences of derived categories 
        and K$3$ surfaces, J. Math.\ Sci.\ (New York) \textbf{84} 
        (1997), 1361-1381. 
\bibitem[Po]{Po} D. Pollack, Algebras and their automorphisms,
        Comm.\ Algebra \textbf{17} (1989), 1843-1866.
\bibitem[RD]{RD} R.\ Hartshorne, ``Residues and Duality,'' Lecture
        Notes in Math.\ {\bf 20}, Springer-Verlag, Berlin, 1966.
\bibitem[Rd]{Rd} J.\ Rickard, Derived equivalences as derived
        functors, J.\ London Math.\ Soc.\ \textbf{43} (1991), 37-48.
\bibitem[Rl]{Rl} C.M. Ringel, ``Tame Algebras and Integral
        Quadratic Forms,'' Lecture Notes in Math.\ \textbf{1099},
        Springer-Verlag, Berlin, 1984.
\bibitem[Rn]{Rn} C.\ Riedtmann, Algebren, Darstellungsk\"{o}cher,
        Ueberlagerungen und zur\"{u}ck, Comment.\ Math.\ Helvetici
        \textbf{55} (1980), 199-224.
\bibitem[RZ]{RZ} R. Rouquier and A. Zimmermann, Picard groups for 
        derived module categories, preprint.
\bibitem[Ye]{Ye} A.\ Yekutieli,
        Dualizing complexes, Morita equivalence and the derived
        Picard group of a ring, 
        J.\ London Math.\ Soc.\ \textbf{60} (1999), 723-746.
\bibitem[YZ]{YZ} A.\ Yekutieli and J.J.\ Zhang, Rings with Auslander
        dualizing complexes, J. Algebra \textbf{213} (1999), 1-51.
\bibitem[Zi]{Zi} A. Zimmermann, 
        Derived equivalences of orders, 
        Canadian Math.\ Soc.\ Conference Proceedings \textbf{18} 
        (1996), 721-749.
\end{thebibliography}
\end{document}